\documentclass[12pt]{article}
\usepackage{graphicx}
\usepackage{dcolumn}
\usepackage{bm}

\usepackage{amsmath}

\newcommand {\hide}[1]        {{ } }

\newcommand{\R}{{\cal{R}}}
\newcommand{\M}{{\cal{M}}}

\newcommand{\I}{{{I}}}

\begin{document}

\bibliographystyle{plain}

\title{Combinatorial Games with a Pass:\\ A dynamical systems approach}

\author{Rebecca E. Morrison\\
Institute for Computational Engineering and Sciences\\
The University of Texas at Austin \\ Austin, TX 78712\\
\it{\{email:  rebeccam@ices.utexas.edu\}} \\
\and Eric J. Friedman\\
School of Operations Research and Information Engineering \\
Cornell University \\ Ithaca, NY \\ \it{\{email:  ejf27@cornell.edu\}}
\and Adam S. Landsberg\\
Keck Science Department \\
Claremont McKenna, Pitzer and Scripps
Colleges \\ Claremont, CA 91711\\ \it{\{email: alandsberg@jsd.claremont.edu\}}}
\date{\today}
\maketitle

\begin{abstract}
By treating combinatorial games as
dynamical systems, we are able to  address a longstanding open question in combinatorial game theory, namely, how the introduction of a ``pass'' move into a game affects its behavior.  We consider two well known combinatorial games, 3-pile Nim and 3-row Chomp.  In the case of Nim, we observe that the
introduction of the pass dramatically alters the game's underlying structure, rendering it considerably more complex, while for Chomp, the pass move is found to have relatively minimal impact. We show how these results can be understood by recasting these games as dynamical systems describable by dynamical recursion relations.  From these recursion relations we are able to identify underlying structural
connections between these ``games with passes'' and a recently introduced class of ``generic (perturbed) games.''  This connection, together with a (non-rigorous) numerical stability analysis, allows one to understand and predict the effect of a pass on a game.
\end{abstract}

\baselineskip = 1.5 \normalbaselineskip

Combinatorial games like Chess, Checkers, Go, Nim, and Chomp  have been the focus of considerable attention in the fields of computer science, mathematics, artificial intelligence, and most recently, chaos and dynamical systems theory.  In traditional combinatorial games (under ``normal play''), two players alternate moves until one player reaches a terminal position
from which no legal move is available, thereupon losing the game \footnote{For a basic introduction to combinatorial games see \cite{WW} or \cite{ANW07}.}. An intriguing but surprisingly difficult question in combinatorial game theory centers on
what happens when standard game rules are modified so as to allow
for a one-time pass -- i.e., a pass move which may be used at most once in a game, and not from a terminal position.
Once the pass has been used by either player,
it is no longer available. Although this question
has been raised in various contexts (see, e.g., \cite{WW, HN03}), it touches upon some deep issues relating to the underlying structure and computational complexity of a game, and
to date it remains largely unanswered.  Indeed, the effect of
a pass on even the simplest, most canonical of combinatorial games -- Nim -- remains an important
open question in combinatorial game theory that has defied traditional approaches, and the late mathematician David Gale even offered a monetary prize to the first person to develop a solution for 3-pile Nim with a pass \cite{Gale}. In this paper we show how tools from dynamical systems theory (wherein we treat ``games with passes'' as a type of dynamical system) can be used to address such issues.

\begin{figure}[h]
\centering
\includegraphics[width = .7\textwidth]{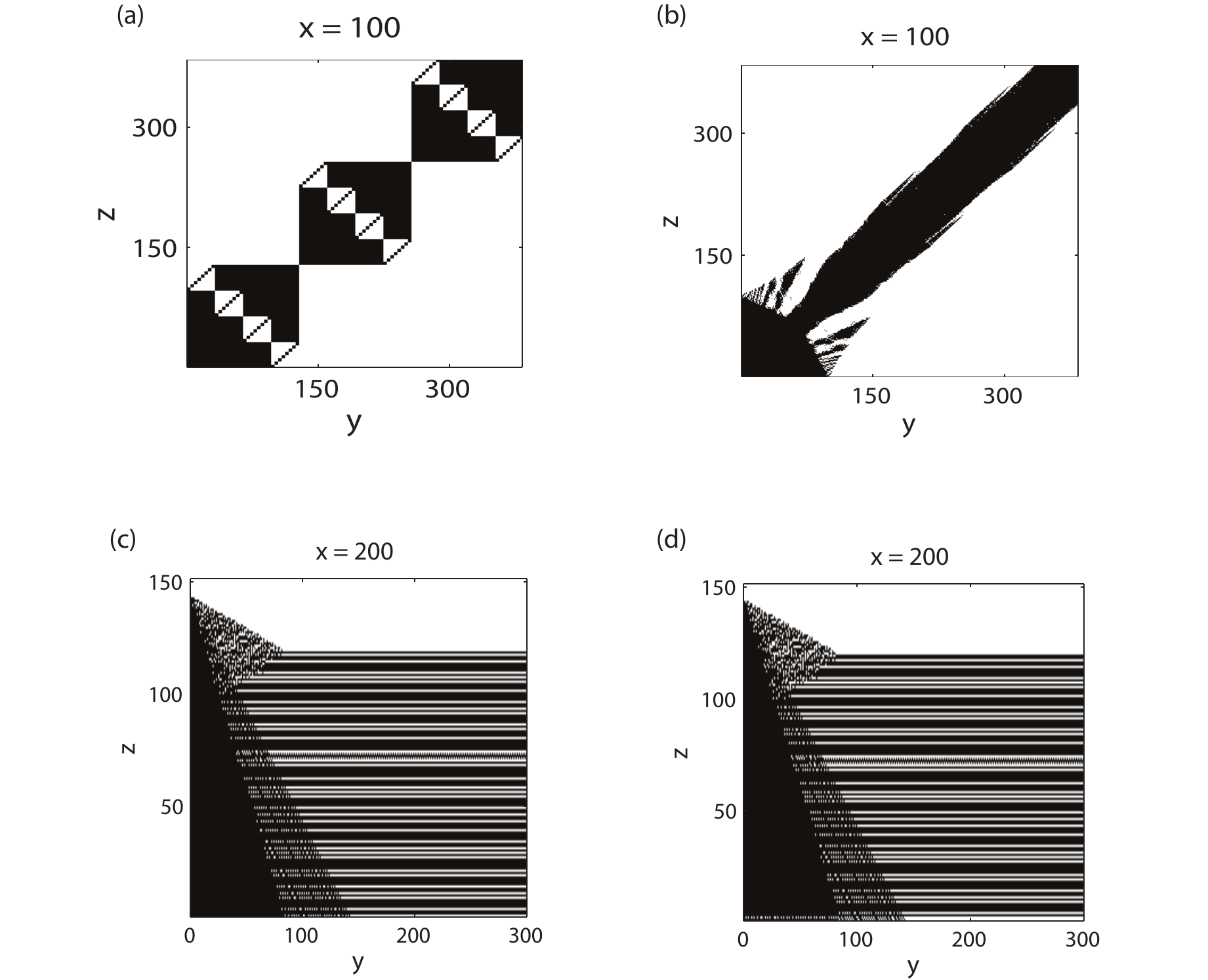}
\caption{The underlying (instant-winner-sheet) geometry of (a) 3-pile Nim, (b) 3-pile Nim-with-a-pass, (c) 3-row Chomp, (d) 3-row Chomp-with-a-pass.}
\label{f1}
\end{figure}

We take up this question of the effects of a pass via two well studied combinatorial games, 3-pile Nim and 3-row Chomp.  The first of these games, 3-pile Nim, is a simple combinatorial game which has been fully solved (without the pass); a complete solution was presented by Bouton over a century ago \cite{Nim}.  The second, 3-row Chomp (without the pass), is an unsolved, intrinsically more complex combinatorial game \cite{Zeil1}.  We find that the introduction of a pass has dramatically different effects on these two games. In the former, the pass leads to a radical change in the game's underlying structure and complexity, while in the latter we find no such dramatic changes (see, e.g., Fig. \ref{f1}, which will be explained later).

 In this paper we attempt to analyze and explain this phenomena
not by offering  an explicit ``solution'' to Nim-with-a-pass or Chomp-with-a-pass -- indeed our analysis indicates why this may not be possible -- but rather
by showing that these combinatorial games with passes are related to a class of ``perturbed" (also called ``generic") combinatorial games \cite{FrL07,FrL08no}.  This transcription
of the problem from one of passes to one of perturbations is based on our discovery of an underlying equivalence in the recursive algorithms used to analyze passes and perturbations in games.   The main focus of this paper is to describe this underlying connection.  This finding will in turn provide
critical insight into why a simple game like Nim is radically transformed by the introduction of a pass, while a game like Chomp is not.  The essential feature of our approach is the recognition that combinatorial games like those considered here can be recast as dynamical systems described by a form of recursive mapping, allowing us to invoke tools and notions from dynamical systems theory such as attractors and sensitivity to initial conditions.

The paper is organized as follows:  In Sect.\ 1  we describe a recently developed recursive technique for analyzing combinatorial games that will also prove useful for analyzing games with passes.  This methodology is based on the construction of recursion operators that describe a game's underlying geometric structure \cite{FrL07, FrL08no}.   In Sect.\ 2 we discuss some general considerations regarding games with passes, and describe how the recursive approach can be adapted to handle passes.  In Sect.\ 3 we analyze the game of 3-pile Nim-with-a-pass.    In Sect.\ 4 we focus on 3-row Chomp-with-a-pass.   In Sect.\ 5  the underlying connections between these two games with passes and so-called ``generic games'' is revealed, and   we discuss the implications of these relationships.

\section{Overview:  Recursion Methodology}
Many of the standard methods of combinatorial game theory have proven ineffective for analyzing the effects of a pass in a game.  In the case of Nim, for instance, Bouton's \cite{Nim} elegant solution method (along with the more general Sprague-Grundy theory which generalizes it \cite{Sprague,Grundy}) breaks down when a pass is introduced and offers very little additional insight. Here, the failure can be attributed to the fact that a pass move renders the game non-decomposable.  The situation is still more difficult in the case of Chomp, a more complex game for which standard solution methods have  failed even in the absence of a pass.

Given this situation, we will utilize a recently introduced dynamical-systems-based methodology for analyzing games that focuses on identifying the underlying geometry of a game and its recursive structure \cite{FrL07,FrL08no,Zeil1,Zeil2,GM11}.  A fundamental virtue of this approach is that it does not rely on a game being solvable or decomposable, and hence can address complex situations, whether it be a simple game which becomes complex by virtue of the introduction of a pass (as in Nim) or a game which is intrinsically complex even in the absence of a pass (as in Chomp). (In a dynamical systems context, this recursive-based approach is somewhat akin to exploiting the fact that a chaotic system is often describable by a simple, underlying iterative mapping.)

In this section we present a basic conceptual overview of the recursive methodology to be used in this paper, restricting first for simplicity to the case of games without the passes, and then describing  the complete methodology which incorporates the modifications necessary to capture the effects of passes.  (We refer the reader to  \cite{FrL07, FrL08no} for a more detailed introduction to the general dynamical-systems-based approach to combinatorial games.)

\subsection{Games without passes: general considerations}
We begin with some general considerations.  To start, we note that in the two games considered here -- 3-pile Nim and 3-row Chomp -- a game position can be represented by a triplet of integers $[x,y,z]$.   In Nim, these integers specify the number of tokens in each pile; in Chomp they are related to the number of tokens in each row. (The specific rules for each game will be described later.)  Hence both games have a three-dimensional ``position space'' \cite{3Dnote}. If we imagine marking the locations in position space of all P-positions \footnote{In keeping with standard terminology from combinatorial game theory, a P-position is a position from which a player will lose if his/her opponent plays optimally; an N-position is one from which a player will win under optimal play.} of the game, they
form a type of three-dimensional geometric object, which we dub the ``P-set''.
The P-set of the game is the main object of interest; knowledge of a game's P-set allows one to define a winning strategy for the game.  Interestingly,
for many combinatorial games, the P-set turns out to be neither a collection of randomly dispersed points in
position space, nor a completely regular geometric object \cite{FrL07, FrL08no}.  Rather, it lies intermediate between the two: It displays
an overall geometric structure, but with local scatter (disorder). (As will be discussed later, the P-set of 3-pile Nim without a pass proves to be exceptional
in this regard; its P-set geometry
is completely regular.)

\begin{figure}[h]
\centering
\includegraphics[width = .4\textwidth]{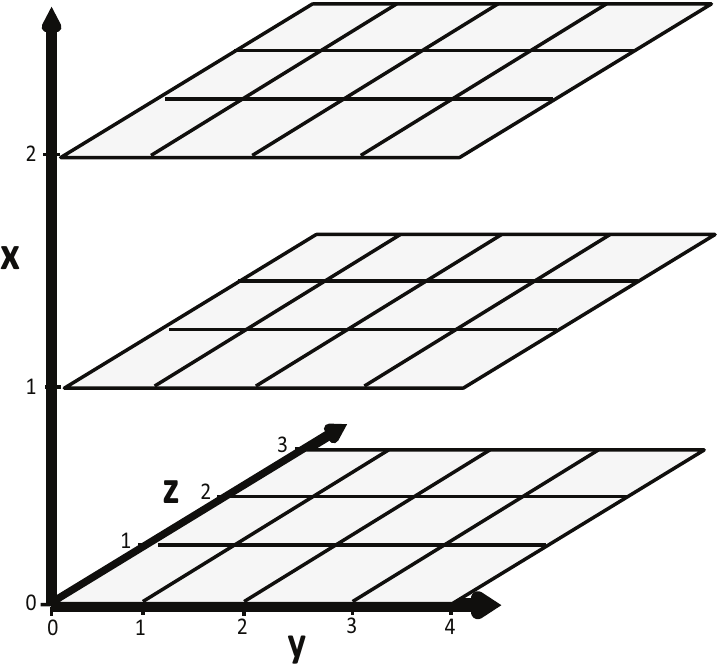}
\caption{Foliating sheets.}
\label{f2}
\end{figure}

A recursion-based analysis of a game's P-set is carried out by foliating the game's three-dimensional position space $x-y-z$ by
sets of  horizontal planes (``sheets''), with $x$ serving as the index for the sheets and $y,z$ as the coordinates within a sheet
(Fig. \ref{f2}).  There are two basic classes of foliating sheets that are of interest: ``loser sheets'' $L_0, L_1, L_2, \ldots$ and ``instant-winner sheets'' $W_1,W_2, W_3, \ldots$.   The loser sheets mark the locations of the game's P-positions.
More formally, the loser sheet at level $x$, denoted $L_x$, is defined to be
a two-dimensional, semi-infinite (boolean) matrix consisting of zeros and ones, with ones
marking the location of the game's P-positions at the specified $x$ value:   $L_x(y,z)=1$ if position $[x,y,z]$ is a P-position
and $0$ otherwise.  We will be interested in the
geometric patterns formed by the P-positions on these loser sheets.
Taken collectively, the loser sheets describe the full geometric structure of the P-set in the three-dimensional position space of the game.

The second class of foliating sheets, the instant-winner sheets, also play a key role in the analysis.
These are defined as follows:  First,
define a game position $[x,y,z]$ to be an {\em instant winner} if
there exists a legal move from that position
to some P-position $[x',y',z']$ on a lower sheet (i.e., with $x'<x$). The instant-winner {\em sheets} are planes in the game's three-dimensional position space on which are marked the
locations of these instant winners. More formally, we define the instant-winner sheet at level $x$, denoted $W_x$, to be
a two-dimensional, semi-infinite (boolean) matrix consisting of zeros and ones, with ones
marking the location of the instant winners having the specified $x$ value:   $W_x(y,z)=1$ if position $[x,y,z]$ is an instant-winning position
and $0$ otherwise. Just as for the loser sheets, we will be interested in geometric patterns that may arise within these instant-winner sheets.

\begin{figure}[h]
\centering
\includegraphics[width = .5\textwidth]{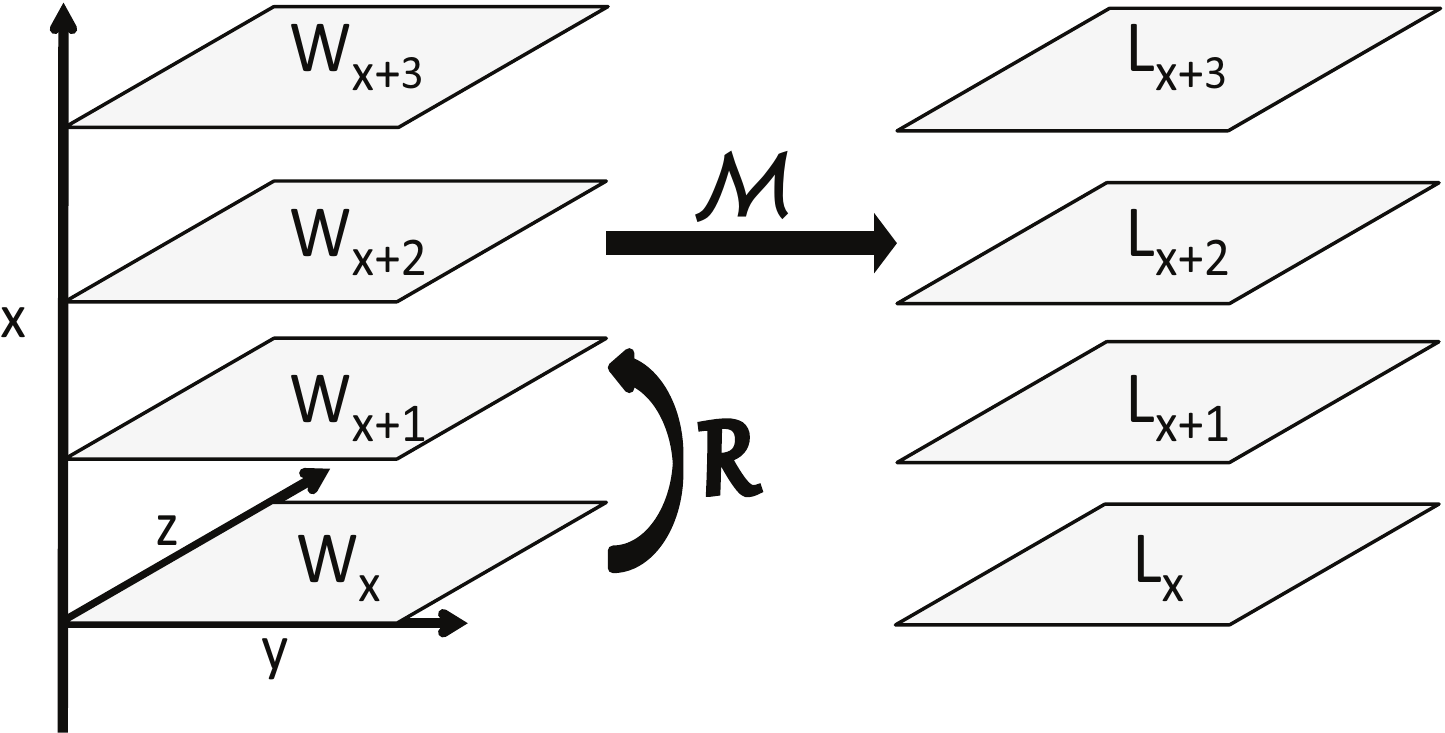}
\caption{Geometric depiction of the action of operators $\R$ and $\M$.    Repeated applications of the recursion operator $\R$ on an instant-winner sheets successively generate the higher-level instant-winner sheets.  Application of the supermex operator $\M$ to a given instant-winner sheet constructs the associated loser sheet.}
\label{f3}
\end{figure}

The basis of the recursive analysis  is the finding that there exist relationships between these various foliating sheets.  (Zeilberger \cite{Zeil1,Zeil2} was among the first to recognize the importance of such recursive structures in games.)  In particular, one can construct a recursion operator, denoted here $\R$, that relates an instant-winner sheet at one level to the instant-winner sheet at the next higher level: $W_{x+1}=\R W_x$. Thus, starting from the lowest-level instant-winner sheet, one can recursively generate all higher instant-winner sheets by repeatedly applying the $\R$ operator.  Moreover, one can  construct another operator, denoted here $\M$ (dubbed the ``supermex operator''), that relates an instant-winner sheet at level $x$ to the associated loser sheet at that same level: $L_x=\M W_x$.  Hence, once the instant-winner sheets have been recursively generated, the loser sheets can  be readily found (and thus the game's P-set determined). This situation is depicted in Fig. \ref{f3}.  For this reason, we can think of the instant-winner sheets as effectively encoding all information about a game's P-set, and thus much of our subsequent analysis will focus on the recursive structure of the instant-winner sheets. We note that the operators $\R$ and $\M$ for the games of 3-pile Nim and 3-row Chomp can be explicitly computed from their respective game rules, as will be illustrated shortly.

\subsection{Games with passes: general considerations}
Adapting the recursive strategy outlined above to games with passes requires some modifications to the operators and foliating sheets, but the overall recursive framework remains intact.  Before discussing this process explicitly, we first mention two basic considerations concerning the introduction of passes into games:

\begin{enumerate}

\item By definition, the pass move cannot be invoked if no tokens remain in the game (i.e.,  more generally, the pass cannot be used from any position that would be deemed a terminal position in the corresponding game without a pass).  Without this restriction, a pass would be uninteresting since a ``game-with-pass" could then be  viewed as  nothing more than a ``game-without-pass" played in disjunctive sum with a nim pile containing a single token -- a situation which is readily handled by standard game-theoretic techniques.

\item Once the pass move has been invoked by either player (thereby rendering it unavailable from that point forward),  the analysis of the game-with-pass reduces to that of the corresponding game-without-pass \footnote{Here, a ``game-without-pass" can refer interchangeably to  either a game (with pass) in which the pass has been invoked by a player and hence is no longer available, or to the original ``pure" game in which a pass move was never allowed in the first place; the two situations are  equivalent.}. Notationally, a hat $\hat{}$ will henceforth indicate that the pass has not yet been used; the absence of a hat signifies that it has been used.  For example, if $\hat{p}$ denotes an arbitrary position in the game when the pass move is still available, $p$ will denote that same position when the pass is no longer available.

\item  If $p$ is a (non-terminal) P-position, then $\hat{p}$ is necessarily an N-position, since by invoking the pass a player at $\hat{p}$ could move his/her opponent to losing position $p$.  Thus, under optimal play, the only time a player at an arbitrary position $\hat{p}$ will invoke the pass move is if position $p$ is a P-position.

\end{enumerate}

The introduction of a pass move into a game also necessitates some modification to the recursive formulation described earlier.   Observe first  that it is possible to continue to view the position spaces of  3-pile Nim-with-a-pass and 3-row Chomp-with-a-pass as three-dimensional   provided we also explicitly specify whether or not the pass remains available.  We do so by using the hat/no-hat notation described above, i.e., if $\hat{p}$ denotes a position with the pass available, $p$ denotes that same position when the pass is unavailable.  In cases where we wish to specify explicitly the coordinates of a position, we do so by appending an auxiliary boolean coordinate as follows: $\hat{p}=[x,y,z;1]$ and $p=[x,y,z;0]$, where the 1 indicates the pass is available and 0 that it has been used.  The introduction of a pass move into a game also requires that we now distinguish between two distinct classes of loser sheets,  $L_0, L_1, L_2, \ldots$ and $\hat{L}_0, \hat{L}_1, \hat{L}_2, \ldots$; the former describe the P-positions of the game at the specified $x$-level when the pass is not available, the latter when the pass is available. A similar distinction must  be made for instant-winner sheets: $W_1, W_2, W_3, \ldots$ and $\hat{W}_1, \hat{W}_2, \hat{W}_3, \ldots$.  More precisely, when the pass is not available, we define a position $p=[x,y,z;0]$  to be an {\em instant winner} if from that position one can move to a P-position $p'=[x',y',z';0]$ at a lower $x$-level (i.e., with $x' < x$). Sheets $W_1, W_2, W_3, \ldots$ mark the locations of such instant winners.   Likewise, when the pass is available, a position $\hat{p}=[x,y,z;1]$ is defined to be an instant winner if one can move to a P-position $\hat{p}'=[x',y',z';1]$ at a lower $x$-level. Sheets $\hat{W}_1, \hat{W}_2, \hat{W}_3, \ldots$ mark the locations of these instant winners in this case.


\section{3-pile Nim-with-a-pass}
In this section we focus on the game of 3-pile Nim-with-a-pass.  We start with a summary of results for the simple case of pure 3-pile Nim (i.e., without the pass), focusing in particular on its underlying geometric and recursive structure.  We then turn to the main focus, Nim-with-a-pass, and analyze its underlying structure.

\subsection{Geometry and recursion for (pure) 3-pile Nim}

In 3-pile Nim, players alternate removing one or more tokens from a pile of their choice, with the player who takes the last token declared the winner.  We let $[x,y,z;0]$ denote an arbitrary position in the game, where $x,y,z$ are non-negative integers representing the number of tokens remaining in each of the three piles. Note here that we are using our extended coordinate representation which explicitly includes a boolean variable (set to be 0) to indicate the absence of a pass.  Although this is redundant in the present context (since pure Nim allows no pass move), it will make for a notationally smoother transition to our subsequent discussion of Nim-with-a-pass.

Given an arbitrary position $[x,y,z;0]$, three types of moves are allowed corresponding to the three choices of piles from which tokens can be removed:
\begin{description}
\item{\textbf{M}1:} \,\, $[x,y,z;0] \longrightarrow [x,y,z-t;0] \,\,\,\,\,0<t\leq z$,
\item{\textbf{M}2:} \,\, $[x,y,z;0] \longrightarrow[x,y-t,z;0] \,\,\,\,\,0<t\leq y$,
\item{\textbf{M}3:} \,\, $[x,y,z;0] \longrightarrow[x-t,y,z;0] \,\,\,\,\, 0<t\leq x.$
\end{description}
We will refer to these accessible positions under these moves as the ``children'' of the original (``parent'') position.  From the game  rules M1-M3 we can construct the recursion and supermex operators for the instant-winner and loser sheets in 3-pile Nim, as we now describe.

To begin, we first introduce  two basic operations on sheets --  the identity operator $\mathcal(I)$ and the addition operator $\mathcal{+}$.  Letting $A$ and $B$ denote two arbitrary sheets (i.e., Boolean matrices), and letting coordinates $(y,z)$ specify an arbitrary element within a sheet, we define:

\begin{description}
\item The $identity$ operator: $\mathcal{I}A=A$
\item The $addition$ operator: $(A\mathcal{+}B)(y,z)=1$ if $A(y,z)=1$ or $B(y,z)=1$.
\end{description}
Note that the addition operator $\mathcal{+}$  is simply the logical OR operation applied to the elements of the sheets.

To compute the instant-winner and loser sheets for 3-pile Nim, recall that an instant-winner sheet $W_x$ consists of all N-positions $[x,y,z;0]$ from which one can move to a P-position at a lower $x$-level, while a loser sheet $L_x$ consists of all P-positions at the specified x value.  From game move M3 above it follows that

\begin{equation}
\label{instantwinner}
W_{x}=L_{0}+L_{1}+\dots+L_{x-1}.
\end{equation}

Next, we construct the ``supermex'' operator $\mathcal{M}$ that relates an instant-winner sheet $W_{x}$ to its associated loser sheet $L_x$ via $L_x=\M  W_x$.  This is done by marching through the elements of the instant-winner sheet one by one (starting from the ``smallest" position under the dictionary ordering and moving toward successively larger positions), locating the first 0 element and marking it as a P-position, then marking off all parents of that P-position on the sheet as N-positions (using games rules M1, M2), then finding the next smallest 0 element, marking it as a P-position, marking off its parents as N-positions, etc.  More formally, the supermex operator $\mathcal{M}$ for the game is defined as the following algorithm:

\begin{description}
\item {\bf Nim supermex operator $\M$ acting on $W_x$}:
{\em
\item 1) Set $\M W_x=0$, $T_x=W_x$, y=0.
\item 2) Let $z_s$ be the smallest $z$ such that $T_x(y,z)=0$ and
set $(\M W_x)(y, z_s)=1$,  $T_x(y+t, z_s )=1$ for all $0\leq t$.
\item 3) Let $y\rightarrow y+1$ and go to step 2.}
\end{description}
By construction,
\begin{equation}
\label{loser}
L_x=\M W_x.
\end{equation}

Finally, simple algebraic manipulations on equations (\ref{instantwinner}) and
(\ref{loser}) yields the desired recursion relation
\begin{equation}
\label{recursion}
W_{x+1}=W_x+\M W_x,
\end{equation}
which, for convenience, can be more compactly expressed as
\begin{equation}
\label{simplerecursion}
W_{x+1}=\R W_x
\end{equation}
with $\R \equiv \I+\M$.

\begin{figure}[h]
\centering
\includegraphics[width = .7\textwidth]{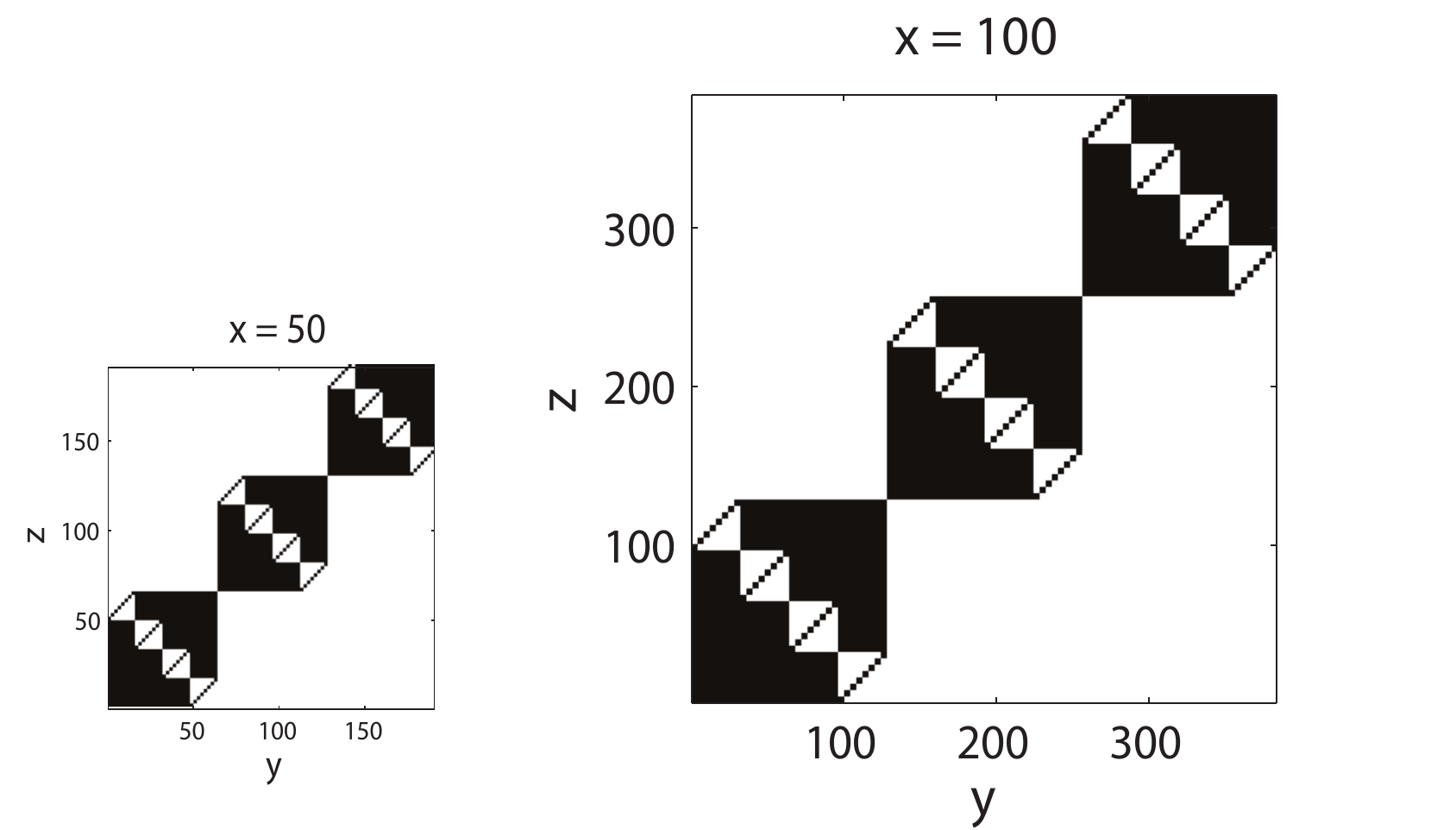}
\caption{Factor-of-two scale invariance for Nim, illustrated with $W_{50}$ and $W_{100}$. (Note that $W_{200}$, $W_{400}, \ldots$, not shown, exhibit this same pattern.)}
\label{f4}
\end{figure}

Recall that the instant-winner sheets $\{W_x\}$ effectively encode all critical information about the game, so we focus attention on them.  A typical instant-winner sheet for pure 3-pile Nim is depicted in Fig. \ref{f1}a.  A key observation is that instant-winner sheets at different $x$-levels have similar overall geometries up to an overall scale factor (which scales linearly with x).  For instance, Fig. \ref{f4} shows instant-winner sheets  $W_x$ for $x=50$ and $x=100$.  The geometric invariance of the sheets is readily apparent.

Three remarks are in order:

\begin{enumerate}
\item Owing to the observed geometric similarity (up to a scale factor) of the sheets, we can conveniently refer to the game's overall geometry without explicitly referencing the particular $x$-level of a sheet.

\item Nim  is a somewhat unusual (in comparison to other games whose instant-winner sheet geometries have been previously investigated) in that the above statement regarding scale invariance of the instant-winner sheets has one subtle complication.  Namely, Nim's instant-winner sheets do not possess a true geometric invariance as suggested above, but rather a type of factor-of-two periodic scale invariance \cite{FrL08no}.  By this we mean that, for {\em any} x,  the  sheets $W_x, W_{2x}, W_{4x}, W_{8x}, \ldots$ will all share the same overall geometry up to rescaling (as illustrated in Fig. \ref{f4} for $x=50$).   Note, however, that two arbitrary sheets $W_x$ and $W_{x'}$, although similar, will not in general have an identical structure (unless $x'$ happens to be a power-of-two multiple of $x$). As an illustration, Fig. \ref{f5} shows instant-winner sheets $W_{x'}, W_{2x'}, \ldots$ for $x'=55$.  The factor-of-two scale invariance is also readily apparent for this set of instant-winner sheets, but note that it is not quite identical to the set shown in Fig. \ref{f4}. For the purposes of the present work, however, the fact that Nim's instant-winner sheets exhibit a {\em factor-of-two} scale invariance rather than a true scale invariance will not prove important for our analysis, and we will continue to refer to Nim's overall ``geometric invariance."

\begin{figure}[h]
\centering
\includegraphics[width = .7\textwidth]{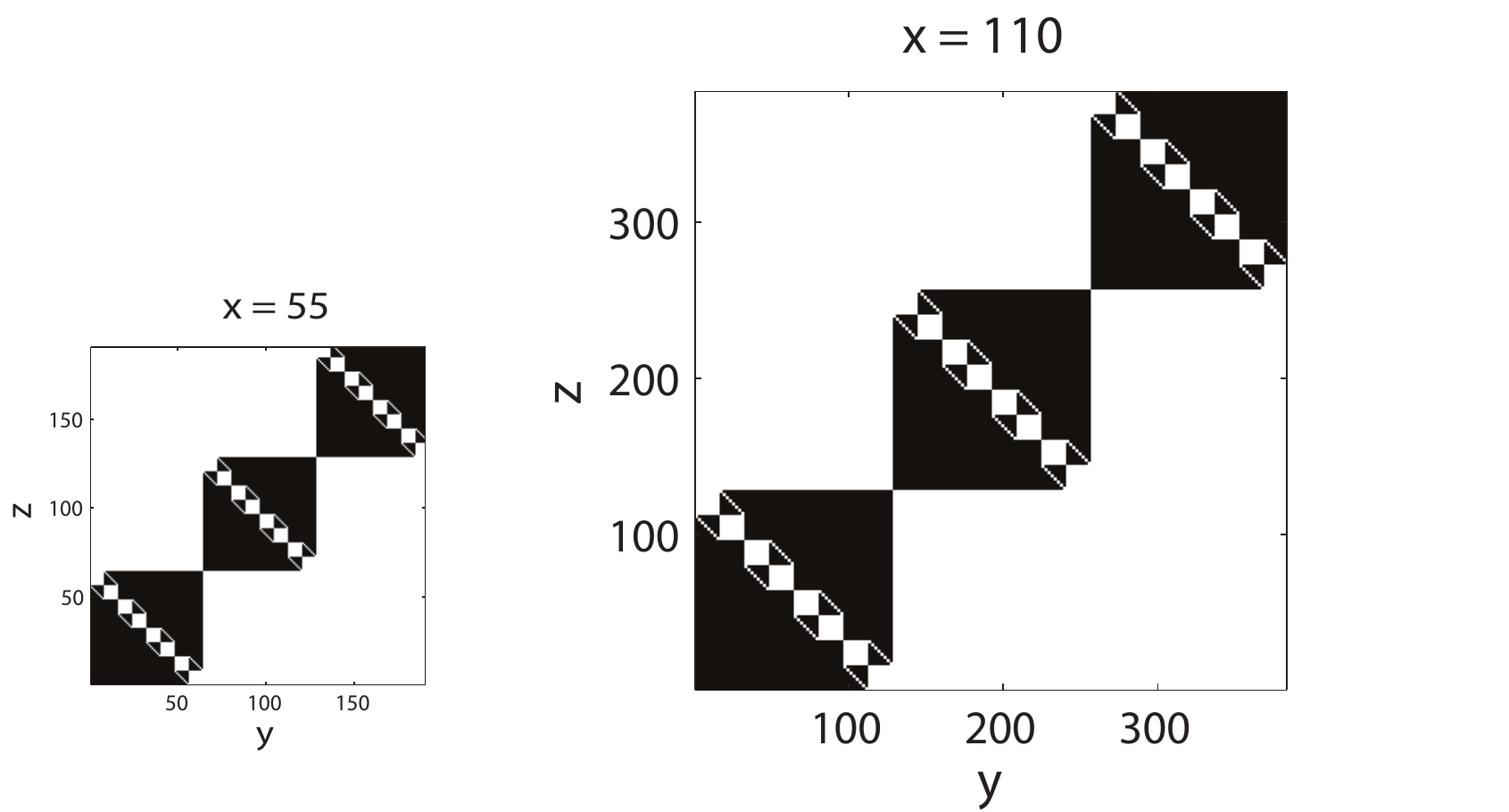}
\caption{Factor-of-two scale invariance for Nim shown, illustrated with $W_{55}$ and $W_{110}$. ($W_{220}$, $W_{440}, \ldots$, not shown, follow this same trend.) Compare with Fig. \ref{f4}.}
\label{f5}
\end{figure}

\item  One of the most prominent and noteworthy features of pure 3-pile Nim (apart from its scale invariance) is the highly regular, ordered nature of its instant-winner-sheet geometry (see, e.g., Figs. \ref{f4} and 5). This is closely associated with the fact that Nim is a solvable game whose P-positions are readily characterized via the theory of nimbers.  In contrast, the geometries of other games we will encounter (i.e., Nim-with-a-pass and Chomp) will prove to be more complex, exhibiting a mixture of ordered features and scatter.
\end{enumerate}

In summary, Eqn.~\ref{recursion} (in conjunction  with Eqn.~\ref{loser}), defines the recursive geometric structure of pure 3-pile Nim. In essence, we have recast this combinatorial game as a dynamical system describable by a type of nonlinear iterative mapping (Eqn.\ \ref{recursion}).

We next turn to the case of Nim-with-a-pass.

\subsection{Geometry and recursion for 3-pile Nim-with-a-pass}

The game rules for 3-pile Nim-with-a-pass are as follows:

\begin{description}
\item{\textbf{M}1:} \,\, $[x,y,z;0] \longrightarrow [x,y,z-t;0] \,\,\,\,\,0<t\leq z$,
\item{\textbf{M}2:} \,\, $[x,y,z;0] \longrightarrow[x,y-t,z;0] \,\,\,\,\,0<t\leq y$,
\item{\textbf{M}3:} \,\, $[x,y,z;0] \longrightarrow[x-t,y,z;0] \,\,\,\,\, 0<t\leq x.$\\

\item{\textbf{M}4:} \,\, $[x,y,z;1] \longrightarrow [x,y,z-t;1] \,\,\,\,\,0<t\leq z$,
\item{\textbf{M}5:} \,\, $[x,y,z;1] \longrightarrow[x,y-t,z;1] \,\,\,\,\,0<t\leq y$,
\item{\textbf{M}6:} \,\, $[x,y,z;1] \longrightarrow[x-t,y,z;1] \,\,\,\,\, 0<t\leq x.$
\item{\textbf{M}7:} \,\, $[x,y,z;1] \longrightarrow[x,y,z;0] \,\,\,\,\,$ provided $[x,y,z] \neq [0,0,0]$.
\end{description}
Here, the last four moves, M4-M7, are associated with the existence of the pass, with M7 in particular signifying the taking of the pass.

\begin{figure}[h]
\centering
\includegraphics[width = .5\textwidth]{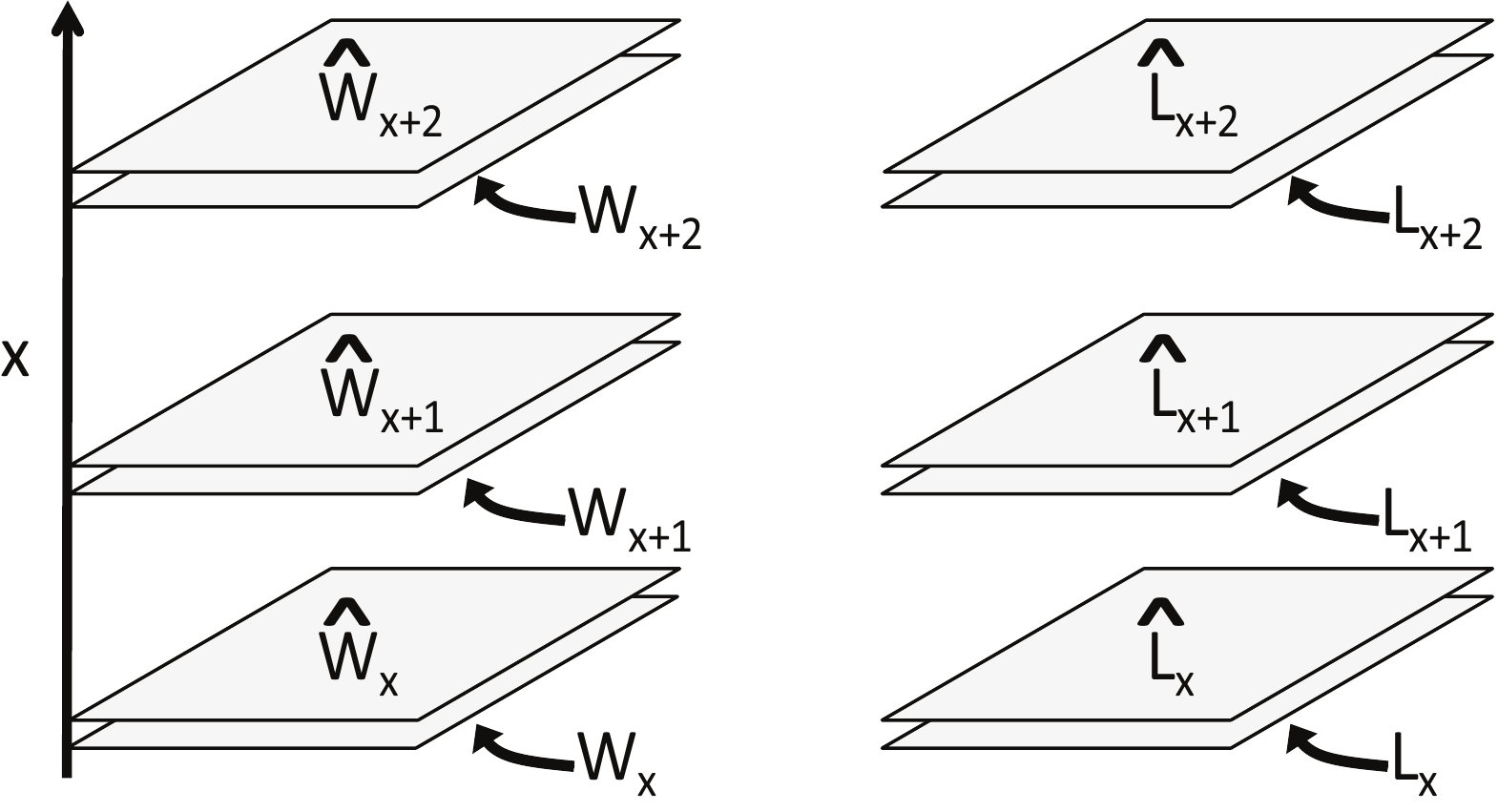}
\caption{The sheet structure for a game-with-pass, analogous to that depicted in Fig. \ref{f3}.  Here, the hatted sheets are for the game-with-pass before the pass has been invoked,  the unhatted sheets for after the pass move has been used.}
\label{f6}
\end{figure}

As previously noted, the introduction of a pass move into a game requires that we now  distinguish between two classes of instant-winner sheets, $W_1, W_2,W_3,\ldots$ and $\hat{W}_1, \hat{W}_2,\hat{W}_3,\ldots$, the former for when the pass has been used, the latter for when it remains available. Similarly we define two classes of loser sheets, $L_0,L_1,L_2,\ldots$ and $\hat{L}_0,\hat{L}_1,\hat{L}_2\ldots$.  Geometrically, it is helpful to visualize these sheets as in Fig. \ref{f6}.  Note that under the game rules  players cannot move from a lower $x$-level sheet to a higher-level sheet, nor can they move from an unhatted sheet to a hatted sheet.

It follows trivially that the structure and sheet operators for the unhatted sheets of Nim-with-a-pass are identical to those of pure Nim.  Namely,
\begin{equation}
\label{nimstructure}
W_{x+1}=W_x+\M W_x \;\;\;  \mbox{and} \;\;\;  L_x=\M W_x.
\end{equation}

The instant-winner sheets for the pass game are also readily constructed.
Recall that an instant-winner sheet $\hat{W}_x$ consists of  all N-positions (at level $x$) with the property that from those positions one can move to  a P-position on a lower-level sheet.   Accordingly, it follows from rule M6 that
\begin{equation}
\label{instantwinnerhat}
\hat{W}_{x}=\hat{L}_{0}+\hat{L}_{1}+\dots+\hat{L}_{x-1},
\end{equation}
in analogy with Eqn.~\ref{instantwinner}.

The construction of the loser sheets differs  from that of pure Nim owing to the introduction of the pass move.   We cannot simply write $\hat{L}_x=\M \hat{W}_x$ (in analogy with Eqn.~\ref{nimstructure}), since this does not incorporate the pass move (rule M7). Specifically, a supermex operator $\M$ acting on solely  the instant winners on sheet $\hat{W}_x$ fails to account for the existence of other N-positions arising directly from the pass move.  These other N-positions, called ``pass-winners,''  are defined as follows:  Position $[x,y,z;1]$ is a {\em pass-winner} iff position $[x,y,z;0]$ is a P-position.  Note that these pass-winners are N-positions created directly by pass rule M7, and are distinct from the instant-winning positions.   Consequently, to construct the loser sheets, we must first add these pass-winners to the instant-winner sheets before applying the supermex operator, as follows:
\begin{equation}
\label{loserhat}
\hat{L}_x=\M (\hat{W}_x + L_x).
\end{equation}
Here, the $L_x$ term represents the contribution from the pass-winners, i.e., the pass-winners at level $x$ are simply given by the P-positions on loser sheet $L_x$.  (To see this, simply note that the pass-winners at level $x$ are  just the parents (under rule M7) of the  P-positions on $L_x$, and hence share the same $[x,y,z]$ coordinates.)  Adding these pass-winners to the instant-winner sheet leads to Eqn.~\ref{loserhat}.

Combining Eqns.~\ref{instantwinnerhat}, \ref{loserhat} yields the desired recursion relation
\begin{equation}
\label{recursionhat}
\hat{W}_{x+1}=\hat{W}_x+\M (\hat{W}_x + L_{x}).
\end{equation}
Together, Eqns.~\ref{recursionhat} and \ref{loserhat} (along with Eqn.~\ref{nimstructure}) define the recursive structure of 3-pile Nim-with-a-pass.  Here, it is most helpful to regard Eqn.~\ref{recursionhat} as an iterative mapping for the instant-winner sheets $\{\hat{W}_x\}$, treating $L_x$ simply as a fixed, known quantity (which can be computed directly from Eqn.~\ref{nimstructure}).

\begin{figure}[h]
\centering
\includegraphics[width = .6\textwidth]{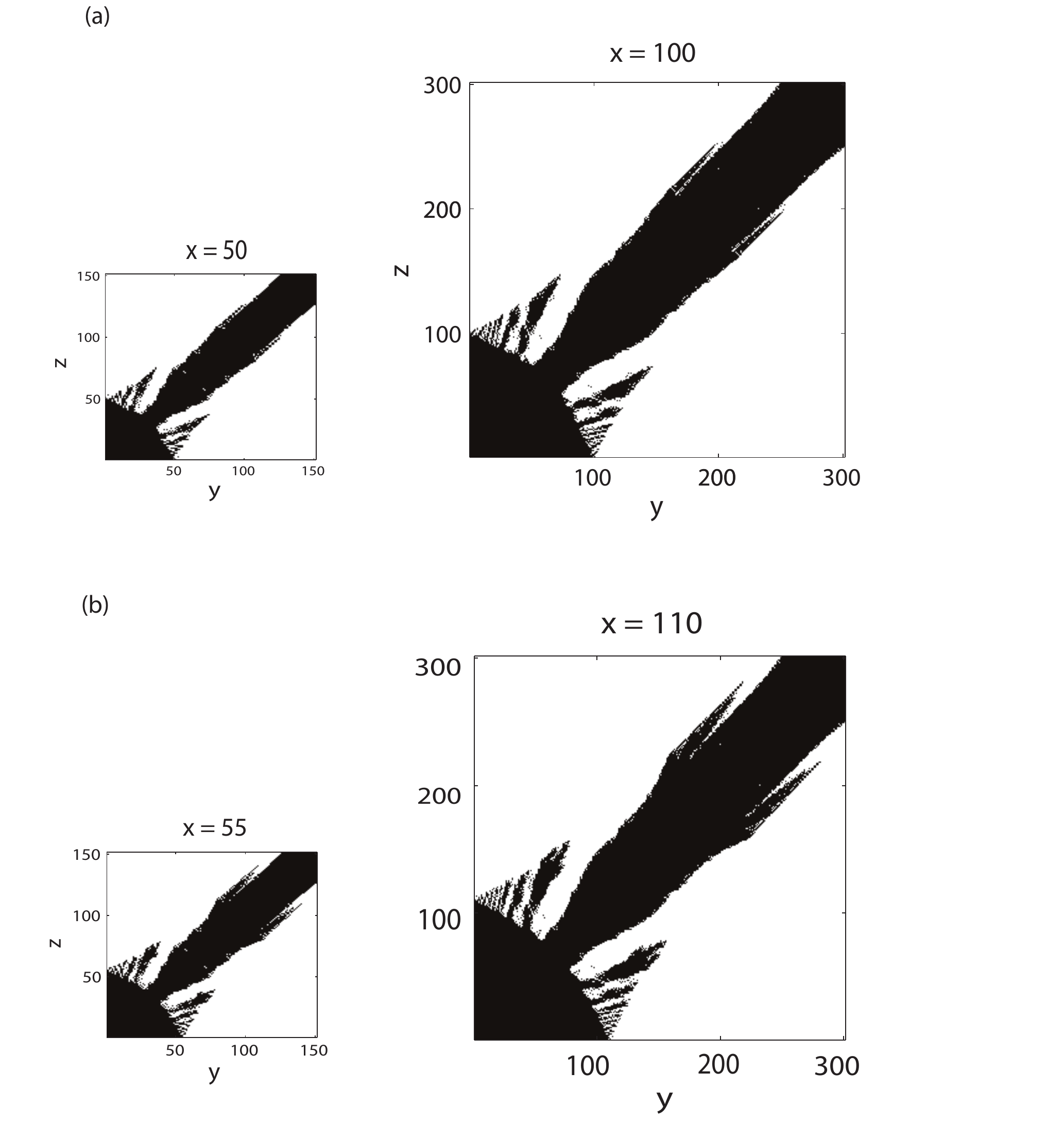}
\caption{The overall instant-winner-sheet geometry for Nim-with-a-pass, shown for (a) $W_{50}$ and $W_{100}$, and (b) $W_{55}$ and $W_{110}$.  Note the overall scale-invariant nature of the geometry;  the slight differences between (a) and (b) reflect the factor-of-two nature of this scale invariance, akin to what was seen in pure Nim. }
\label{f7}
\end{figure}

Fig. \ref{f1}b depicts the instant-winner-sheet geometry for Nim-with-a-pass, as generated from the above recursion relations. Several features are worthy of note.  First, we find that the instant-winner sheets of Nim-with-a-pass display a remarkable geometrical invariance -- i.e., although there is not a precise, one-to-one correspondence between points on different sheets, the
{\em overall} geometry of sheets at different $x$-levels is the same up to a scale factor (which increases linearly with $x$), as illustrated in Fig. \ref{f7}a. (We mention though that, just as for pure Nim, it is not a true geometric invariance in that there exists an underlying factor-of-two periodicity, as can be seen in a comparison of Figs. \ref{f7}a and 7b; however, this subtlety will not prove important for the present analysis.)  Second, note the striking difference in overall geometry between pure Nim and Nim-with-a-pass.  The introduction of the pass renders the geometry considerably more complex, destroying the regularity seen in the case of pure Nim.

The dramatic change in geometric structure caused by the pass has many implications that will be discussed in more detail later.  From the outset, though, it leads to an important inference which will guide our overall approach to the problem.  The issue centers on the likelihood of uncovering an analytical solution to the game of 3-pile Nim-with-a-pass.  In particular, we note that to date no analytical solution to the game (analogous to Bouton's solution for pure Nim) has ever been found, and whether this is even possible has remained an open question.  One might reasonably conjecture that perhaps this failure arises because  the introduction of a pass increases the game's computational complexity in a fundamental way, rendering it analytically intractible.   If so, we would naturally expect that such an increase in complexity would be manifest in the game's underlying instant-winner-sheet geometry.  Reasoning conversely, given the geometrically complex appearance of the game's instant-winner sheets observed in Fig. \ref{f7}, which stands in marked  contrast to pure Nim's highly regular geometry, we are led to conjecture that perhaps an analytical  ``solution" (\`{a}
la Bouton) to Nim-with-a-pass does not exist, or, at a minimum, would be  inherently difficult to discover.   Accordingly, in this paper we will not
attempt to find an analytical solution to the problem of Nim-with-a-pass, but will instead pursue an alternate approach to understanding the effects of a pass by relating the game's recursive structure to that of another class of games.   Before doing so, however, we will first show how the general recursive approach leading to Eqn.~\ref{recursionhat} for Nim-with-a-pass can be readily applied to other games with passes as well.  For this purpose we consider next the game of 3-row Chomp-with-a-pass and work out its recursive structure.  Doing so will yield new insights into the  effects of passes in games.  We will then show how in fact both Nim-with-a-pass and Chomp-with-a-pass are intimately related to a class of so-called ``generic games'' and can be understood in this context.

\section{3-row Chomp-with-a-pass}
In this section we focus on the game of 3-row Chomp-with-a-pass.  We start with a review of some known results for the case of pure 3-row Chomp (i.e., without the pass), focusing in particular on its underlying geometric and recursive structure.  We then turn to the main focus of this section, Chomp-with-a-pass.

\subsection{Geometry and recursion for (pure) 3-row Chomp}

\begin{figure}[h]
\centering
\includegraphics[width = .6\textwidth]{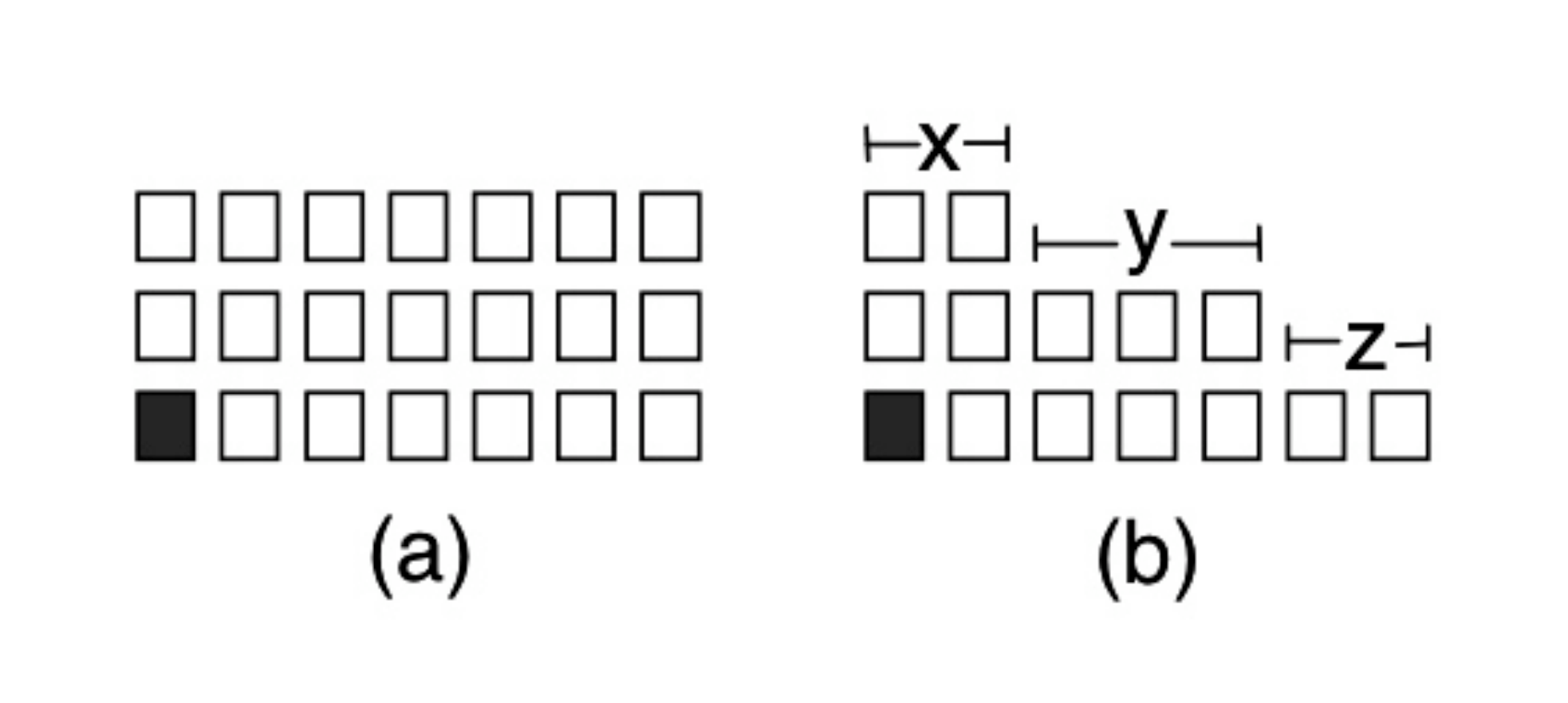}
\caption{3-row Chomp.  (a) The starting configuration.  The darkened token in the lower-left corner represents the poison token. (b) An intermediate position, described in general by coordinates $[x,y,z]$.}
\label{f8}
\end{figure}

3-row Chomp represents one of the simplest ``unsolved'' combinatorial games.  It is played with three rows of tokens (Fig. \ref{f8}a).  On his/her turn, a player selects a token and removes that token along with all tokens above or to the right of the selected token. Play alternates between the two players until eventually one player is forced to take the ``poison'' token located in the lower left corner, thereby losing the game.

At any stage of play, the current position of the game can be conveniently described by a triplet of integers $[x,y,z]$, as illustrated in Fig. \ref{f8}b \cite{Zeil1}.  Although we consider here only the pure game of Chomp (i.e., without a pass), we will adopt the extended notation  $[x,y,z;0]$, with the Boolean coordinate set to 0 to indicate the absence of a pass.

The allowed moves in Chomp are as follows:

 \begin{description}
 \label{purechompmoves}
\item{\textbf{M}1:} \,\, $[x,y,z;0] \longrightarrow [x,y,z-t;0] \,\,\,\,\,0<t\leq z$,
\item{\textbf{M}2:} \,\, $[x,y,z;0] \longrightarrow[x,y-t,z+t;0] \,\,\,\,\,0<t\leq y$,
\item{\textbf{M}3:} \,\, $[x,y,z;0] \longrightarrow[x,y-t,0;0] \,\,\,\,\, 0<t\leq y.$
\item{\textbf{M}4:} \,\, $[x,y,z;0] \longrightarrow [x-t,y+t,z;0] \,\,\,\,\,0<t\leq x$,
\item{\textbf{M}5:} \,\, $[x,y,z;0] \longrightarrow[x-t,0,z+y+t;0] \,\,\,\,\,0<t\leq x$,
\item{\textbf{M}6:} \,\, $[x,y,z;0] \longrightarrow[x-t,0,0;0] \,\,\,\,\, 0<t\leq x.$\end{description}
Position $[0,0,1;0]$, wherein only the poison token remains, is treated as a terminal position.

Construction of a recursion operator for Chomp proceeds analogously to what was done for Nim, though the detailed process is somewhat more complicated.  To begin, we first define some general operations on sheets which will prove useful.  In addition to the identity operator and Boolean addition operator defined previously for Nim, we must also introduce a left-shift operator $\cal{L}$ and diagonal-element-adding operator $\cal{D}$ as follows.  Letting $A$ denote an arbitrary sheet (i.e., Boolean matrix) and $(y,z)$ the coordinates within a sheet, define:
\begin{description}
\item The {\em left-shift} operator $\cal{L}$: $\mathcal{L}A(y,z)=A(y+1,z)$
\item The {\em diagonal-element-adding} operator $\cal{D}$: If $A$ is a loser sheet and letting $A(0,z^*)$ denote the (unique) non-zero element of $A$ at $y=0$, define
$\mathcal{D} A(y,z)=A(y,z)+\delta_{z,z^*-y}$, where $\delta$ denotes the Kronecker delta function.
\end{description}
An instant-winner sheet $W_x$ can be expressed as
\begin{equation}
\label{chompIN}
W_x=\mathcal{L}^x \mathcal{D} L_0+\mathcal{L}^{x-1} \mathcal{D} L_1 + \ldots + \mathcal{L} \mathcal{D} L_{x-1}.
\end{equation}
The terms in this expansion are analogous to those in Eqn.~\ref{instantwinner} for Nim, i.e.,  the first term represents the contribution to $W_x$ from the parents of P-positions on loser sheet $x=0$, the second term from parents of P-positions on loser sheet $x=1$, etc.  These  follow directly from Chomp game rules M4 and M5 together with the definition of an instant-winner sheet; see \cite{FrL07} for a more detailed derivation.  For present purposes though, it suffices to simply note that each term in Eqn.~\ref{chompIN} represents a contribution from a loser sheet at a particular $x$-level.

Just as for Nim, we can construct a supermex operator for Chomp that relates an instant-winner sheet at level $x$ to the corresponding loser sheet:
\begin{equation}
\label{chomploser}
L_{x}=\mathcal{M}W_{x}.
\end{equation}
Note here that the structure of this equation is identical to Eqn.~\ref{loser} for Nim, though the definition of the supermex operator $M$ differs for the two cases since the game rules are different.  In Chomp, the supermex operator follows directly from game rules M1-M3 above. Explicit details can be found in \cite{FrL07}; for present purposes, however, the general form of Eqn.~\ref{chomploser} will suffice.

Combining the above two relations yields the desired recursion relation for Chomp:
\begin{equation}
\label{chomprecursion}
W_{x+1}=\mathcal{L}W_x+\mathcal{L}\mathcal{D}\mathcal{M}W_x,
\end{equation}
or, equivalently,
$W_{x+1}=\mathcal{R}W_{x}$, where $\mathcal{R}\equiv \mathcal{L}(\mathcal{I}+\mathcal{D} \mathcal{M}).$

A representative instant-winner sheet $W_x$ is shown in Fig. \ref{f1}c. Just as in Nim, these sheets exhibit a remarkable geometric scale invariance -- i.e., all
sheets in Chomp share a common overall geometric pattern (whose size increases linearly with $x$). As before, note that {\em geometrical invariance} of the sheets refers to their overall structure -- i.e., to the fact that their slopes of boundary lines and density of points in the various regions are all the same -- not that there exists an exact, point-by-point correspondence between points on different sheets.  Observe also that Chomp's instant-winner-sheet geometry exhibits both regular features (e.g., distinctive regions with well-defined linear boundaries) and disorder (e.g., scattered points within regions). We also remark that, unlike for Nim,  a true geometric invariance for Chomp holds
for all sheets -- i.e., the subtle powers-of-two periodicity noted for Nim is absent.
A full discussion of the geometry of pure Chomp can be found in \cite{FrL07}.

\subsection{Geometry and recursion for 3-row Chomp-with-a-pass}
We now consider the introduction of a pass move into Chomp, and seek to determine a recursive formulation of this problem.  We begin with the game rules.  Here, the original rules M1-M6 of pure Chomp must be supplemented with additional rules associated with the pass:

 \begin{description}
 \label{passchompmoves}
\item{\textbf{M}7:} \,\, $[x,y,z;1] \longrightarrow [x,y,z-t;1] \,\,\,\,\,0<t\leq z$,
\item{\textbf{M}8:} \,\, $[x,y,z;1] \longrightarrow[x,y-t,z+t;1] \,\,\,\,\,0<t\leq y$,
\item{\textbf{M}9:} \,\, $[x,y,z;1] \longrightarrow[x,y-t,0;1] \,\,\,\,\, 0<t\leq y.$
\item{\textbf{M}10:} \,\, $[x,y,z;1] \longrightarrow [x-t,y+t,z;1] \,\,\,\,\,0<t\leq x$,
\item{\textbf{M}11:} \,\, $[x,y,z;1] \longrightarrow[x-t,0,z+y+t;1] \,\,\,\,\,0<t\leq x$,
\item{\textbf{M}12:} \,\, $[x,y,z;1] \longrightarrow[x-t,0,0;1] \,\,\,\,\, 0<t\leq x.$
\item{\textbf{M}13:} \,\, $[x,y,z;1] \longrightarrow [x,y,z;0] \,\,\,\,\,$ provided $[x,y,z] \neq [0,0,1]$.
\end{description}

Once the pass has been invoked (rule M13), Chomp-with-a-pass reduces to pure Chomp, and hence the previous recursion relations (Eqns.~\ref{chomploser}, \ref{chomprecursion}) apply.  While the pass is still available, it follows directly from rules M7-M9 that the instant-winner sheets may be written as
\begin{equation}
\label{chomppassIN}
\hat{W}_x=\mathcal{L}^x \mathcal{D} \hat{L}_0+\mathcal{L}^{x-1} \mathcal{D} \hat{L}_1 + \ldots + \mathcal{L} \mathcal{D} \hat{L}_{x-1},
\end{equation}
in analogy with Eqn.~\ref{instantwinnerhat} for the pure game.
However, the construction of the loser sheets will differ from that of pure Chomp owing to the introduction of the pass move, i.e., the relation $\hat{L}_x=\M \hat{W}_x$ (the analog of Eqn.~\ref{chomploser}) does not hold since it does not incorporate the pass move (rule M13). In particular, one must add the pass-winners to the instant-winner sheet before applying the supermex operator.   (Recall that a position $[x,y,z;1]$ is a {\em pass-winner} iff position $[x,y,z;0]$ is a P-position; these pass-winners are N-positions created directly by pass rule M13.) Thus we have
\begin{equation}
\label{chomppassloser}
\hat{L}_x=\M (\hat{W}_x + L_x),
\end{equation}
with $L_x$ term representing the contribution from the pass-winners.

Combining the above expression with Eqn.~\ref{chomppassIN} yields the desired recursion relation
\begin{equation}
\label{chomppassrecursion}
\hat{W}_{x+1}=\mathcal{L} \hat{W}_x+\mathcal{L} \mathcal{D} \mathcal{M}(\hat{W}_x + L_{x}).
\end{equation}

Eqns.~\ref{chomppassrecursion}, \ref{chomppassloser}, in conjunction with \ref{chomploser}, \ref{chomprecursion}, constitute the basic recursion relations for the game of 3-row Chomp-with-a-pass.  Here, it is most helpful to regard Eqn.~\ref{chomppassrecursion} as an iterative mapping for the instant-winner sheets $\{\hat{W}_x\}$, treating $L_x$ simply as a fixed, known quantity (which can be computed directly from Eqn.~\ref{chomploser}).

A typical instant-winner-sheet geometry, computed from these recursion relations, is shown in Fig. \ref{f1}d.  Most strikingly, observe that the overall geometric structure of 3-row Chomp-with-a-pass appears to be {\em identical} to that of pure 3-row Chomp \footnote{Here as before, {\em identical} means that the various slopes and densities in different sheets are the same, not that there is an exact, point-by-point correspondence.}. Note that this similarity arises despite the differences in the underlying recursion relations for Chomp-with-a-pass and pure Chomp.  This stands in marked contrast to what was observed previously for the cases of pure Nim and Nim-with-a-pass, wherein the introduction of the pass gave rise to a radically different instant-winner-sheet geometry.

A preliminary insight into why the cases of Nim and Chomp might differ in their response to a pass is related to the solvability of the game:  Recall that, unlike for pure Nim, there is no known analytical solution of pure Chomp, suggesting more  generally that perhaps the introduction of a pass move into solvable games tends to significantly modify their underlying geometry (i.e., analytical tractability is a delicate, structurally unstable feature in games), whereas more complex games which lack an analytical solution are intrinsically more robust and will not be radically modified by introduction of a pass.
Though a potentially interesting line of pursuit and closely related to previous work \cite{FrL08no}, in this paper we take a different approach to understanding the effects of passes, and will demonstrate the existence of deeper connections between games with passes and so-called generic games, as we describe next.

\section{Generic games}

\subsection{Introduction to generic games}
The notion of a ``generic'' (aka ``perturbed'') game was introduced in \cite{FrL07}. The idea is straightforward.  Given a combinatorial game,  rather than considering this game in isolation, instead regard it as one member of a larger family of games (the``generic class'' of the game) formed by making various rule changes to the original game.  Specifically, consider a rule change of the following type:  Select an arbitrary P-position in the original game and simply declare it to be an automatic N-position -- i.e., a player landing on this chosen position is declared the winner in this (modified) game.
In effect, we have created a variant of the original game through this simple rule modification. More generally, one could select other P-positions for conversion, or even multiple P-positions simultaneously (though for now we restrict the number of such modifications to be finite). Related to this is the analysis in \cite{FO98} which considered three specific sets of such perturbations in a generalized Wythoff's game. Different choices of  P-positions selected for conversion will lead to different variants of the original game.  Collectively, we will refer to the set of all variants of the original game as the ``generic'' class of the game. (The reason for this terminology will become clear momentarily.)

A complementary view of this  notion of generic games comes from a dynamical-systems theory perspective:  In dynamical systems theory, one often considers an iterative mapping $T: \Re^n \rightarrow \Re^n$ taking a point $x \in \Re^n$ to a new point  $x \rightarrow Tx.$  Here, one  focuses on the orbit \{$x_0, Tx_0, T^2x_0, T^3x_0,\ldots$\}
of  an initial point $x_0$ under repeated iterations of this mapping.  Of particular interest are the asymptotic behavior of this orbit and the degree to which this asymptotic behavior is affected by choice of initial condition $x_0$. We can think of changing an initial condition $x_0$ as essentially ``perturbing" the system. If  the asymptotic behavior of an orbit of a dynamical system is stable to changes in initial conditions then the system is said to possess an ``attractor'' (see, e.g., \cite{GuH86}).  Certain attractors are known to exhibit a type of ``sensitivity to initial conditions'' under such perturbations and are thus deemed {\em chaotic}, a notion which can be put on a quantitative footing by computation of the Lyapunov exponents associated with the orbit.

By recasting a combinatorial game in terms of a recursion relation,  it becomes possible to treat the game as a type of dynamical system.  In particular, a recursion operator on the instant-winner sheets (e.g., Eqn.~\ref{simplerecursion}) can be regarded as an iterative mapping that relates an instant-winner sheet at one level to the instant-winner sheet at the next level up. In this context, the ``mapping" acts on instant-winner sheets (i.e., semi-infinite Boolean matrices) rather than points in $\Re^n$, but is qualitatively similar.  Moreover, we can now think of the rule changes that gave rise to the variant games in the generic class as corresponding to perturbations on the instant-winner sheets.  More precisely, declaring a P-position at level $x$ to be an automatic N-position creates a perturbation of the next instant-winner sheet $W_{x+1}$ (see Eqn.~\ref{simplerecursion}).
Thus, each variant game is associated with a particular perturbation of the instant-winner sheets, and is the analogue of perturbing the initial condition $x_0$ in a dynamical system.  So the study of the generic class of a game reduces to an analysis of a dynamical system $(W \rightarrow \R W)$ under changes to initial conditions.  This in turn allows one to import dynamical-systems notions like attractors, sensitivity to initial conditions and chaos into the realm of generic combinatorial games, as described in \cite{FrL07}.

Before exploring this notion of a generic class of a game in analytical and computational detail, we first present a key observation about generic games that suggests a possible link between these generic games and games-with-passes.

\subsection{An observation about generic games}

\begin{figure}[h]
\centering
\includegraphics[width = .6\textwidth]{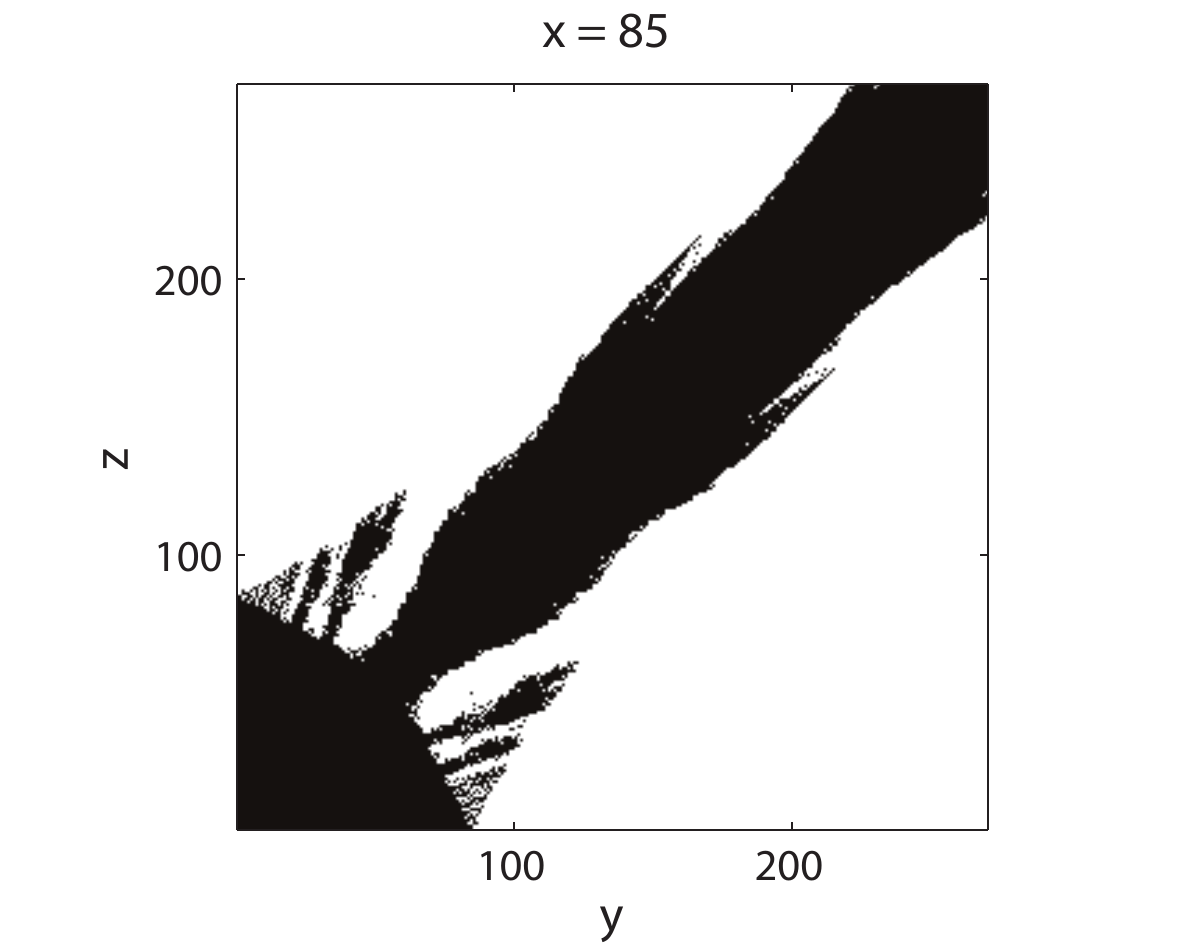}
\caption{Geometry of Generic Nim, illustrated for $W_x$ with $x=85$.  As noted, other choices of $x$-values exhibit this same basic geometry (to within the factor-of-two scale invariance).}
\label{f9}
\end{figure}

The first hint of a potential connection between the notion of the generic class of a game described above and a game-with-a-pass comes from a consideration of Nim. We can  create an instantiation of Generic Nim by making a simple rule change to pure Nim -- e.g., select an arbitrary P-position, say $[9,13,4]$, and  declare it to be an automatic N-position.
We can then recursively construct the instant-winner sheets for this variant game.  Fig. \ref{f9} shows the typical instant-winner-sheet geometry that results (illustrated with $x=85$).

Several preliminary remarks are in  order.

\begin{enumerate}
\item The overall geometry seen in the figure is representative of all $x$-levels of this instancing of Generic Nim, not just $x=85$ (save for the subtle factor-of-two nature of the scale invariance noted previously).

\item Had we created a different instantiation of Generic Nim by making some other arbitrary rule change (e.g., declaring, say, P-position $[7,4,3]$ to be an automatic N-position), we would find that the {\em same} overall geometry seen in Fig. \ref{f9} again emerges.   In fact, numerical simulations suggest that most variants of Nim (i.e., nearly all instantiations of Generic Nim) share the same overall geometry seen in the figure (up to a factor-of-two periodicity) despite the fact that each variant is generated via a different rule change.  This remarkable commonality is the underlying reason for our use of the terminology {\em generic} to describe this class of variant games.  We note that although this observation is based on numerical simulations for randomly generated rule changes (i.e., random perturbations of the initial conditions) and has not been formally proven, it nonetheless does appear that the generic class of the game comprises {\em a single large attractor which contains most variants of the game}.

Alternatively, rather than emphasizing  the limiting behavior of the instant-winner sheets themselves, it is also useful to consider the limiting behavior of the  ``smoothed'' probability distribution of the instant-winner sheets, as discussed in \cite{FrL07,FrL08no}. By ignoring the point-by-point details and focusing on the probabilistic structure of the sheets we see that the limiting behavior of the probability distribution, in this sense, is that of an attractive fixed point.  In other words, our numerical experiments suggest that  ``smoothed asymptotics'' of Generic Nim is a single stable fixed point \footnote{There are some additional subtleties here, due to the infinite-dimensional nature of the problem (as will be discussed later) and the existence of zero measure neutral directions.}.

\item An analogous situation holds for Generic Chomp as well.  In particular, we observe that virtually all instantiations of Generic Chomp appear to share the same underlying geometry as one another.  Alternatively stated, the numerical experiments suggest that  ``smoothed asymptotics'' of Generic Chomp is a single stable fixed point, i.e., the variant games comprising Generic Chomp look identical.  The geometry of these variant games is shown in Fig. \ref{f10}.

\begin{figure}[h]
\centering
\includegraphics[width = .6\textwidth]{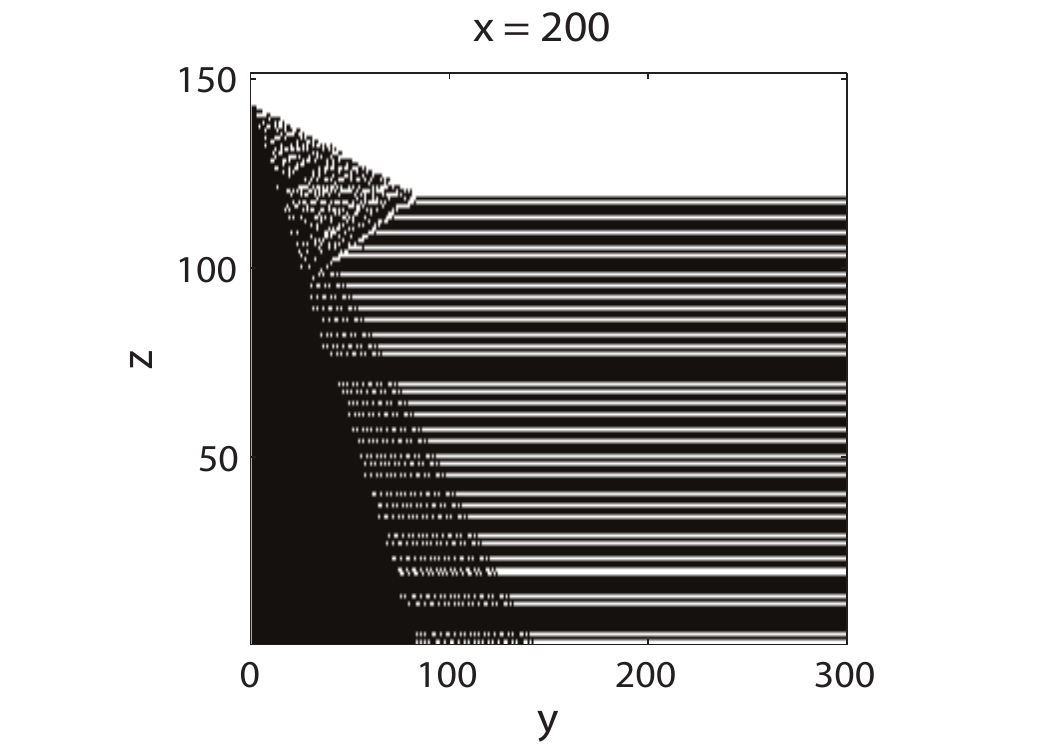}
\caption{Geometry of Generic Chomp, illustrated with $W_{200}$.}
\label{f10}
\end{figure}


\end{enumerate}

For the present purposes, however,  the most striking and important observation stemming from Fig. \ref{f9} is the remarkable similarity between the instant-winner-sheet geometry of Generic Nim (Fig. \ref{f9}) and that of Nim-with-a-pass (Fig. \ref{f1}b).  Likewise, we find that Generic Chomp (Fig. \ref{f10}) and Chomp-with-a-pass (Fig. \ref{f1}d) also share the same geometry as one another \footnote{Quantitatively, the similarity of the sheets can be assessed in several ways. In the case of Chomp, for example,  the overall sheet geometry is characterized by six fundamental geometric parameters (associated with the properties of various boundary lines), as described in \cite{FrL07}.  We have numerically computed these parameters  for Generic Chomp and Chomp-with-a-pass, and have found them be the same to within numerical uncertainty.}. These tantalizing geometric similarities  in fact provide the first hint of a possible deep structural connection between  generic games and games-with-passes,  and indeed this numerical observation is what originally motivated the present study.  Understanding this connection is the primary focus of the remainder of this paper.

To put this relationship on an analytical footing, we will derive recursion relations for  Generic Nim and Generic Chomp.  A  direct comparison with the corresponding recursion relations for Nim-with-a-pass and Chomp-with-a-pass will unveil the sought-after connection.  Before proceeding, however, we  first mention an important subtlety which will arise.  We have defined generic games via a {\em finite} number of rule changes.  However, in the following we will need to consider a generalization  which allows for an infinite number of perturbations. The distinction between the finite and infinite cases may appear superficial at first blush  -- and indeed the recursive formulation for the two cases is identical --  but some care is nonetheless required.  In particular,  while a finite number of perturbations can be naturally associated with changes to the instant-winner sheets (i.e., initial conditions), an infinite number of perturbations will turn out to be most naturally viewed as a change to the recursion operator $\R$ itself.  We will to use the term {\em generic}  to refer to either case, but will use {\em finite}-generic or {\em infinite}-generic when it is necessary to explicitly distinguish between the two cases.

\subsection{Recursive formulation of Generic 3-pile Nim}
In this section we derive recursion operators for Generic (3-pile) Nim. As noted, the rules of Generic Nim are identical to those of pure Nim, save for the  modification that certain specified P-positions in the pure game are now to declared to be automatic N-positions --  i.e., the game terminates anytime such a position is reached, with the player at that position declared the winner.  In what follows we  consider the most general case wherein we allow for the possibility that multiple P-positions on either a finite or infinite number of sheets have been selected for conversion to automatic (terminal) N-positions.

Towards this end, we begin by introducing a set of ``variant sheets'' $\{V_x\}, x=1,2,3,\ldots$ on which are specified the locations of the P-positions (of pure Nim) that have  been declared to be automatic  N-positions of the variant game.  More precisely, the variant sheet at level $x$, denoted $V_x$, is defined to be a two-dimensional, semi-infinite (boolean) matrix consisting of zeros and ones, with ones marking the location of the designated positions:  $V_x(y,z) =1$ if position $[x,y,z]$ has been designated for conversion to an automatic N-position and 0 otherwise.

Next, we define instant-winner and loser sheets for Generic Nim just as was done for pure Nim.  Notationally, we will use  tildes,
$\{\tilde{W}_{x}\}, \{\tilde{L}_{x}\},$  when referring to the  sheets of the generic game, so as  to distinguish them from those of pure Nim and Nim-with-a-pass.   The definitions of these sheets are entirely analogous to earlier definitions: A loser sheet at level $x$ of Generic Nim, denoted $\tilde{L}_x$, is a two-dimensional semi-infinite boolean matrix such that $\tilde{L}(y,z)=1$ if $[x,y,z]$ is an P-position and 0 otherwise. An instant-winner sheet at level $x$ of Generic Nim, denoted $\tilde{W}_x$, is a two-dimensional semi-infinite boolean matrix such that $\tilde{W}(y,z)=1$ iff $[x,y,z]$ is an N-position from which one can move to a P-position at a lower $x$-level.   Finally, note that the rules of Generic Nim are identical to rules M1-M3 of pure Nim, save for the added stipulation that if position $[x,y,z]$ is one of the designated automatic (terminal) N-positions (i.e., $V_x(y,z)=1$), then no moves from this position are allowed.

We next construct the instant-winner sheets for Generic Nim.  Just as in pure Nim, a sheet $\tilde{W}_x$ contains all the
 parents of the P-positions at the lower $x$-levels. Thus we have:
\begin{equation}
\label{gnimIN}
\tilde{W}_{x}=\tilde{L}_{0}+\tilde{L}_{1}+\dots+\tilde{L}_{x-1}.
\end{equation}

The loser sheets of Generic Nim are constructed in a manner similar to those of pure Nim via the standard supermex operator $\M$, but must appropriately account for not only the instant winners but also the designated automatic N-positions.  In particular, we must add the variant sheet to the instant-winner sheet before applying the supermex operator:
\begin{equation}
\label{gnimloser}
\tilde{L}_{x}=\mathcal{M}(\tilde{W}_{x}+V_x).
\end{equation}
Combining the preceding two expressions yields the desired recursion relation for Generic Nim:
\begin{equation}
\tilde{W}_{x+1}=\tilde{W}_{x}+\mathcal{M}(\tilde{W}_{x}+V_{x}). \label{gnimINrecursion}
\end{equation}
Hence, given sheet $\tilde{W}_x$ at level $x$, one can generate the sheet at the next higher level. (Here, we are treating the variant sheet $V_{x}$ as a fixed, known quantity, since it contains the designated points which define the variant game.) Thus, from repeated applications of Eqns.~\ref{gnimloser}, \ref{gnimINrecursion}, we can recursively determine all P-positions of Generic Nim.

Comparing these recursion relations for Generic Nim with those derived earlier for Nim-with-a-pass, a striking structural similarity emerges:
\begin{align}
\text{{\em Generic Nim:}}  \;\;\;\tilde{W}_{x+1}&=\tilde{W}_{x}+\mathcal{M}(\tilde{W}_{x}+V_{x}), \nonumber \\
\text{{\em Nim-with-a-pass:}} \;\;\;\tilde{W}_{x+1}&=\tilde{W_{x}}+\mathcal{M}(\tilde{W}_{x}+L_{x}). \nonumber
\end{align}
We see in fact that, in terms of their underlying dynamical recursion relations, Nim-with-a-pass represents a specific instantiation of Infinite-generic Nim, wherein the variant sheets $\{V_{x}\}$ of Generic Nim are selected to be precisely the loser sheets $\{L_{x}\}$ of Nim-with-a-pass.
Before discussing the implications of this central result, we first show that an analogous finding holds for the game of Chomp, namely, Chomp-with-a-pass proves to be a specific instantiation of Infinite-generic Chomp.

\subsection{Recursive formulation of Generic 3-row Chomp}
The construction of a recursive formulation for Generic Chomp proceeds similarly to that of the preceding case.  In analogy with Eqn.~\ref{chompIN} for pure Chomp, the instant-winner sheets can be expressed as
\begin{equation}
\label{gchompIN}
\tilde{W}_x=\mathcal{L}^x \mathcal{D} \tilde{L}_0+\mathcal{L}^{x-1} \mathcal{D} \tilde{L}_1 + \ldots + \mathcal{L} \mathcal{D} \tilde{L}_{x-1}.
\end{equation}
 (Here again we employ the tilde notation $\tilde{W}_x, \tilde{L}_x$ to signify explicitly  that we are working with the generic class of the game.)  Owing to the presence of the designated automatic N-positions, the loser sheets of Generic Chomp are given by
\begin{equation}
\tilde{L}_x=\M(\tilde{W}_x+V_x), \nonumber
\end{equation}
where the supermex operator $\M$ is identical to that of pure Chomp.
Combining the preceding two equations yields
\begin{equation}
\label{gchompresult}
\mbox{{\em Generic Chomp:}} \;\; \; \; \tilde{W}_{x+1}=\mathcal{L}\tilde{W}_{x}+{\mathcal{L}\mathcal{D}\mathcal{M}(\tilde{W}_{x}}+V_{x}).
\end{equation}
Comparison with our earlier result
\begin{equation}
\label{chomppassrecursionagain}
\mbox{{\em Chomp-with-a-pass:}} \;\;\;\; \hat{W}_{x+1}=\mathcal{L} \hat{W}_x+\mathcal{L} \mathcal{D} \mathcal{M}(\hat{W}_x + L_{x})
\end{equation}
reveals that Chomp-with-a-pass is in fact a specific instantiation of Infinite-generic Chomp, wherein the variant sheets of Infinite-generic Chomp $\{V_{x}\}$ have been selected so as to coincide with the loser sheets   of Chomp-with-a-pass $\{L_{x}\}$.

\section{Analysis and Discussion}
The recursion-operator formulation of games-with-passes has revealed that they share an underlying structural connection with  infinite-generic games, culminating in the relations defined by Eqns.~\ref{recursionhat}, \ref{gnimINrecursion}, and \ref{gchompresult}, \ref{chomppassrecursionagain}.  This constitutes one of the central findings of this paper.

With this discovery of an analytical connection between these two classes of games in hand, one would like to extend this result and claim  that a game-with-a-pass can in fact be fully understood as a specific instantiation of an infinite-generic game (with the study of passes thereby effectively subsumed by the more general study of generic games).  However, some caution is warranted here, as several issues must first be addressed before this can be done.  As is well known from dynamical systems theory, merely because the structure of the dynamical recursion operators for one system can be expressed as a subcase of the other (e.g., compare Eqns.~\ref{gchompresult} and \ref{chomppassrecursionagain}, with the assignment $V_{x} \rightarrow L_{x}$) does not imply that the dynamical behavior of the subcase is necessarily  identical to that of the generic case, and thus invoking the general properties of the generic case to ``explain'' the properties of the pass case could fail.  Such an occurrence  is commonly seen in the study of dynamical systems.  Consider, for instance, the simple iterative mapping $x \rightarrow 2x$.  Though the orbit of virtually any initial value $x_0$  will run off to $\pm \infty $ under repeated iterations of the mapping, there is an exceptional subcase, namely $x_0=0$, for which the generic behavior does not hold.  Likewise, the structural similarity of the recursion operators for games-with-passes and infinite-generic games does not guarantee that the former will necessarily inherit the asymptotic dynamical characteristics of the latter, i.e., there is no {\em a priori} guarantee that the game-with-pass subcase will lie in the basin of attraction of the attractor associated with the generic game.

Given these and related considerations, several numerical issues must be addressed, including
\begin{enumerate}
\item To what extent is the generic class of Nim or Chomp truly comprised of a single large attractor (i.e., could the set of initial conditions outside of the attractor's basin of attraction have non-zero measure)?
\item How can we quantify whether the instant-winner-sheet geometry of the games-with-passes in fact belongs in the generic class?
\item What is the relationship between an infinite-generic game and a finite-generic game?  In particular, do an infinite number of rule changes lead to different structural properties than a finite number of such changes?  (The relevance of this question will become clear shortly.)
\end{enumerate}

To address these and related issues, we have conducted a variety of numerical experiments.  As will be described next, these include
\begin{enumerate}
\item Simulations of generic Nim and generic Chomp in which large numbers of random perturbations are made to every instant-winner sheet, demonstrating the robustness of their attractors
\item A numerical analysis of the degree of perturbation associated with the pass games, demonstrating that these perturbations are, in a certain sense, small
\item A numerical study of the games' sensitivity to initial conditions under the dynamics of the recursion operators, demonstrating why dynamical correlations between different sheets can be expected to be small

\end{enumerate}

\section{Supporting numerical findings}
\subsection{Robustness to random perturbations}
\begin{figure}[h]
\centering
\includegraphics[width = .6\textwidth]{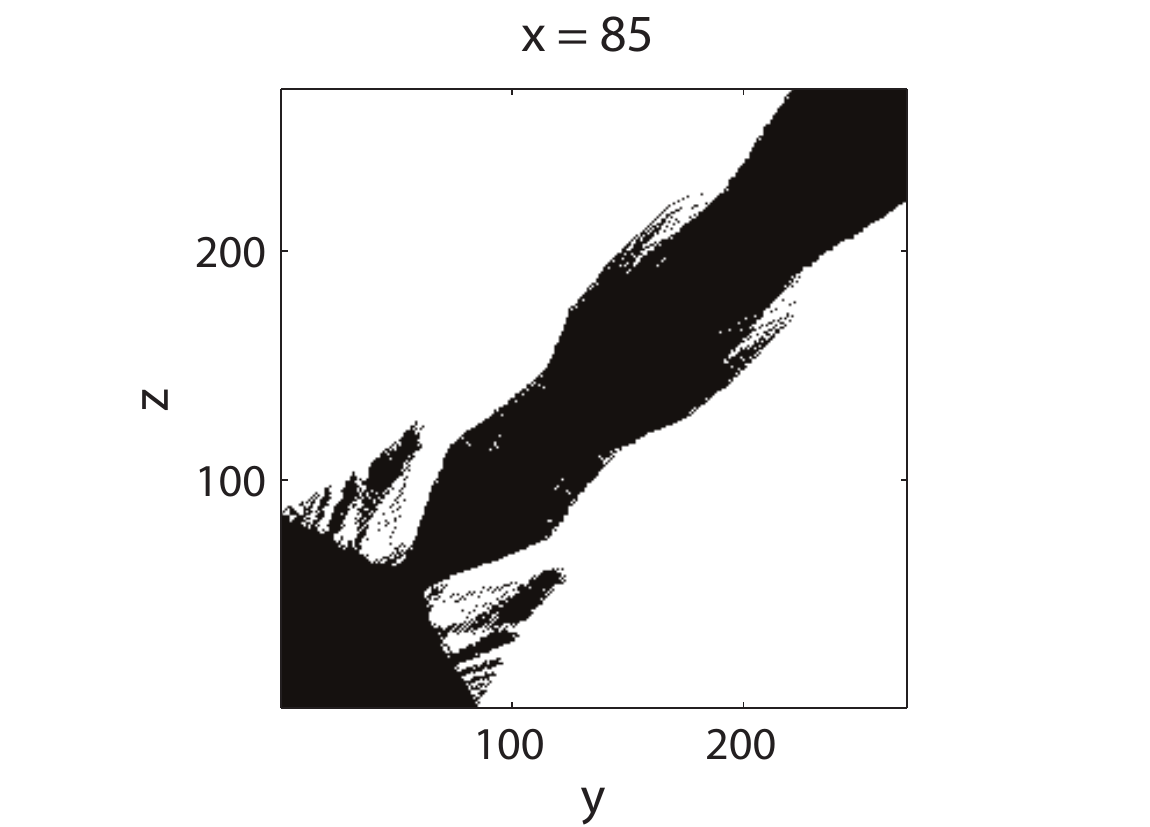}
\caption{The geometry of Infinite-generic Nim, illustrated with $W_{85}$.  The perturbations made to each column were drawn from a normal distribution centered along the main diagonal of the figure.  (We note that the rather modest differences between the geometry of finite-generic Nim (Fig. \ref{f9}) and Fig. \ref{f11} can be largely attributed to the factor-of-two invariance, not the finite versus infinite nature of the perturbations.)}
\label{f11}
\end{figure}

As noted previously, numerical simulations show that if a single perturbation is made at a randomly chosen point on a single instant-winner sheet of either pure Nim or pure Chomp, thereby creating an instantiation of the generic game, then the same overall instant-winner geometry almost always results (e.g., Figs. \ref{f9}, 10) regardless of which initial point was selected as the perturbation site \cite{FrL07}.  This suggests there is a single large attractor for this dynamical system, or, at a minimum, an attractor with a very large basin.  If instead of perturbing at a single point, we instead make perturbations to a finite number of points on a finite number of sheets, we again numerically observe that the same stable attractor re-emerges after iterating the mapping (i.e., recursion relation).  Hence the attractors seen in the generic games seem to be highly robust.  Although it is difficult to accurately quantify the actual rapidity of the approach to the attractor following a perturbation (since calculations of associated eigenvalues or Lyapunov exponents for such games remains an open problem), nonetheless the decay back onto the attractor  is found to occur after only a small number of iterations, and is not inconsistent with an exponential decay \cite{FrL07}.  However, we must now address the issue of an infinite number of perturbations.  In particular, comparing the recursion relations for the games-with-passes with the generic (perturbed) games (Eqns.~~\ref{recursionhat}, \ref{gnimINrecursion}, and \ref{gchompresult}, \ref{chomppassrecursionagain}), we see that the addition of the pass is associated with perturbations {\em at every $x$-level} -- i.e., an infinite number of perturbations in total. Viewed as a dynamical system then, what we have is system with an attractor in which we are making repeated perturbations to the mapping, never allowing the system time to fully recover from a perturbation before hitting it with the next member of this infinite progression of perturbations.  So it is not obvious {\em a prior} whether under these circumstances the original attractor would in fact remain stable.
In addition, we need to also consider the possibility that the attractor itself might change, i.e., that the fixed point of the smoothed problem is changed by the repeated perturbations. To test this, we ran numerical simulations in which a random perturbation was made on {\em every} column of {\em every} instant-winner sheet.  Here, by ``column" we mean at every $y$ value of a sheet (recall that $(y,z)$ are the coordinates on a sheet).  The results for the game of Nim are shown in Fig. \ref{f11}.  Again we find that the same overall geometry re-emerges (up to differences associated with factor-of-two scale invariance).

To understand how this is possible, it is useful to have an analogy in mind.  Consider the one-dimensional diffeomorphism $x \longrightarrow x/2$, which possesses an attracting fixed point at $x=0$.  If one makes a single small perturbation away from the attractor (which can be readily realized by setting the initial condition $x_0$ to a small nonzero value $\epsilon$), then it is clear that the system will exponentially decay back onto the attractor $x=0$.  To mimic the effects of repeatedly perturbing the system away from the attractor (once per iteration), one simply modifies the mapping to $x \longrightarrow x/2 + \epsilon$.  A quick check reveals that even under these repeated perturbations the system remains close to the original attractor (within $2 \epsilon$ in this case).  Hence something akin to this appears to be occurring in the games-with-passes -- the repeated perturbations ($L_x$) to the system at every $x$-level
(see Eqn.~\ref{recursionhat}) caused by the pass move are numerically observed to be insufficient to move the system far from the attractor associated with the generic game. (This is consistent with small perturbations of a stable fixed point.) However, we also find evidence that something still stronger helps keep the system close to the attractor, namely, the effective ``size" of the perturbations themselves decreases with $x$-level, which allows the original fixed point to be unmodified, as we describe next.

\subsection{Degree of perturbation in the pass games}
Comparing Eqns.~\ref{recursionhat}, \ref{gnimINrecursion}, and \ref{gchompresult}, \ref{chomppassrecursionagain}, we have seen that the pass move effectively acts as a perturbation to the system (under the association $V_x \rightarrow L_x$).  So at every $x$-level the system is being modified by the addition of the P-positions of the original game (i.e., $L_x$), which are now being designated to be automatic N-positions.

Focusing on these P-positions for a moment, note that on any $n \times n$ loser sheet $L_x$ there are ${\cal O}(n)$ P-positions.  This follows from the fact that for a given $x$ and $y$,  there is a unique $z$-value for which position $[x,y,z]$ is a P-position. (Suppose, for the sake of contradiction, that $[x,y,z_{1}]$ and $[x,y,z_{2}]$ were both P-positions, with $z_{1}>z_{2}$. Then a player could directly move from one P-position  to another under game rule M1, in violation of a basic principle of combinatorial game theory.)

\begin{figure}[h]
\centering
\includegraphics[width = .7\textwidth]{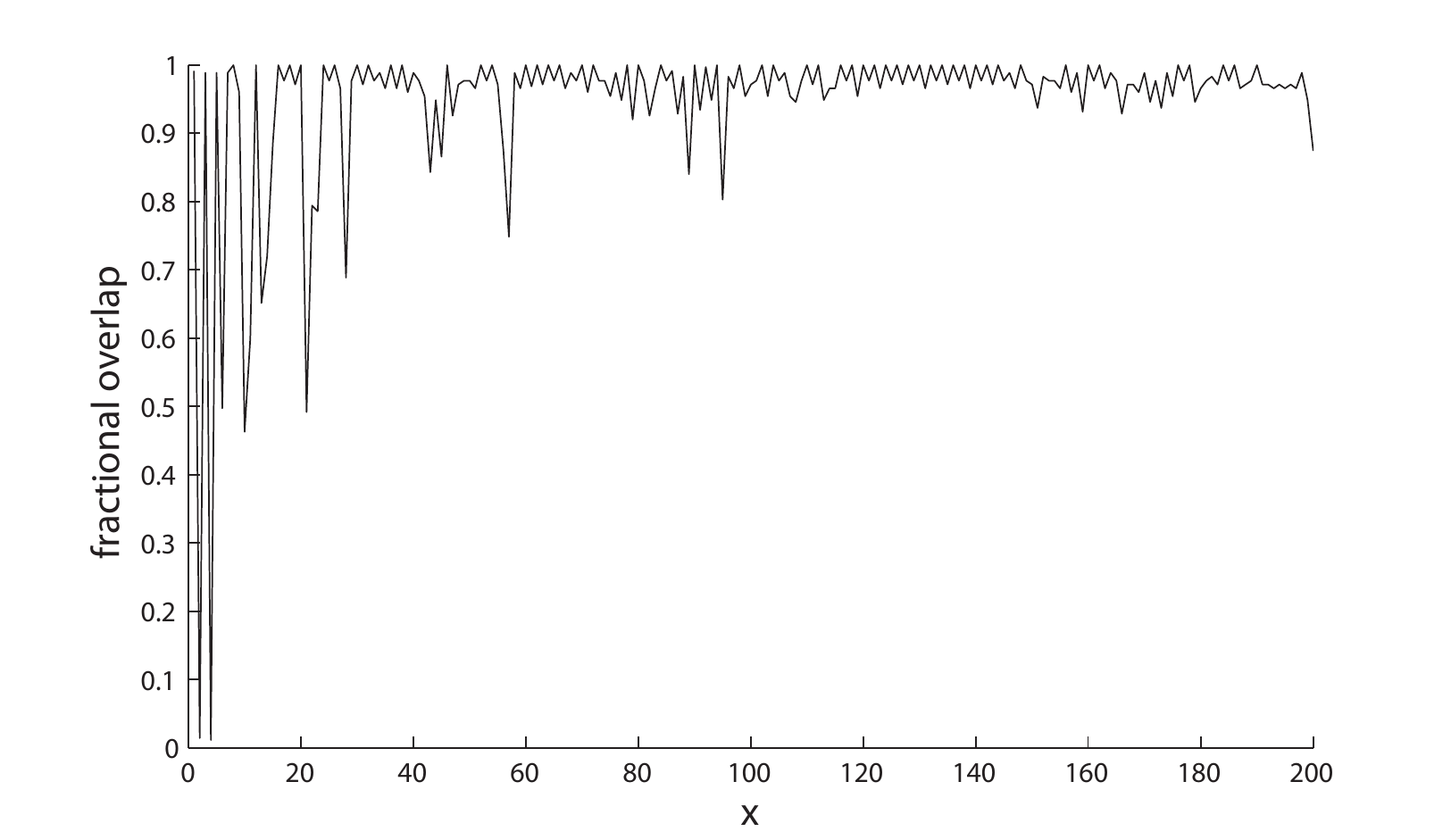}
\caption{At each $x$-level, the plot shows the fractional overlap between the actual P-positions of Nim-with-a-pass ($M(\hat{W}_x+L_x)$) and what the P-positions would have been had the pass-winners $L_x$ at that $x$-level not been included. At low $x$-levels the degree of overlap is small, indicating that the effect of the pass-winners is most pronounced there.  However, the overall upward trend of the plot indicates that the effective strength (on average) of the perturbations introduced by the pass-winners becomes progressively weaker with higher $x$-levels.}
\label{f12}
\end{figure}

So on each $n\times n$ sheet the pass move effectively introduces  ${\cal O}(n)$ perturbations. However, not all of these perturbations actually affect the dynamical behavior of the recursion relations.  To see this, consider creating one instantiation of a generic game in which only a single position (on a single sheet) is designated an automatic N-position.  Note that this solitary perturbation will only have an effect if the position selected for conversion had been a P-position of the game prior to the perturbation (since clearly declaring an existing N-position to be an automatic N-position has no meaningful effect).  Thus, of the ${\cal O}(n)$ individual perturbations on a sheet introduced by $L_x$, only some of these will have actual dynamical consequences, depending upon whether or not they would have been P-positions had the perturbation not been made.  However, determining which of the points in $L_x$ have an effect is somewhat subtle, in that the effect (or lack thereof) of a particular point depends on which other perturbing points are assumed to be present.  To circumvent this subtlety we instead  measure the net influence of the entire set of perturbations $L_x$ by numerically computing the overlap between sheets $M(\hat{W}_x+L_x)$ and $M \hat{W}_x$, where the first term describes actual P-positions of the perturbed (i.e., pass)  game at level $x$, while the second describes what the P-positions would have been had the entire perturbation at level $x$, namely $L_x$, not been made.  In Fig. \ref{f12} we plot this overlap as a function of $x$-level. The key numerical finding here is that the overall influence of the perturbations decreases as a function of x; in actual effect, successive perturbations are becoming progressively weaker.  Hence the attractor seen in the generic game remains stable even in the pass game in part because the effective perturbations introduced by pass move ($L_x$) become progressively weaker as one goes to higher-level sheets. (This analysis is related to that in \cite{FO98} for perturbations in a generalized Wythoff's game.)

\subsection{Sensitivity to initial conditions}
\begin{figure}[h]
\centering
\includegraphics[width = .7\textwidth]{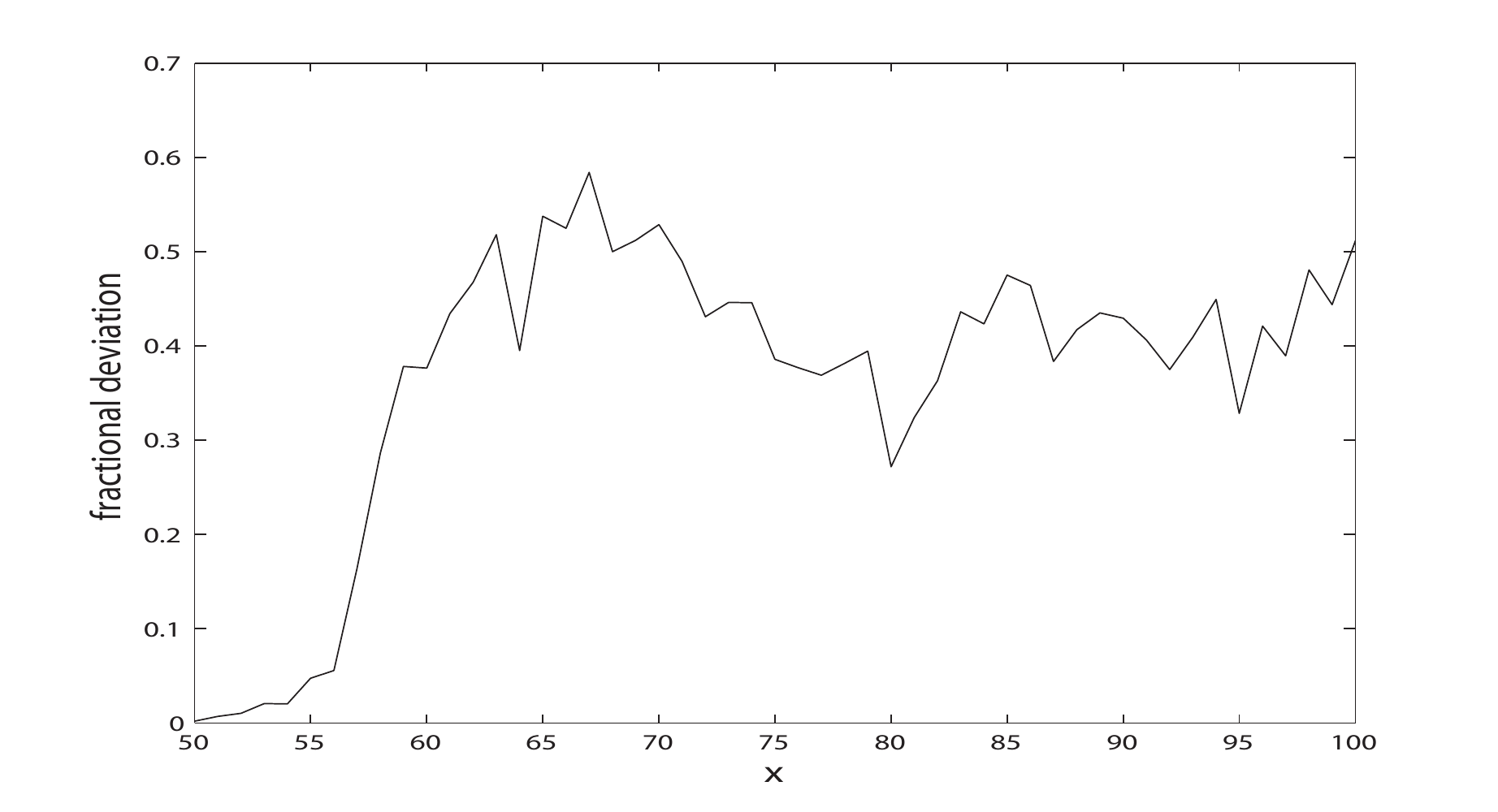}
\caption{Sensitivity to initial conditions in Generic Nim.  The plot was constructed by comparing one instantiation of Generic Nim with a second instantiation created by making a single perturbation to instant-winner sheet $W_{50}$ of the first. The resulting differences in the locations of the P-positions between the two instantiations were then tracked as a function of $x$-level.  The plot shows that these (fractional) differences rise very rapidly from zero -- i.e., after moving up only about 15 levels in $x$, nearly half of all P-positions have changed locations as a result of that single perturbation to the sheet $x=50$.}
\label{f13}
\end{figure}

Underlying the preceding discussion concerning why the pass game exhibits the same geometric attractor seen in the generic game is a presumption that the perturbations introduced by pass move are not so highly correlated with one another that they could  yield a different geometric structure than that of the generic game. Partial numerical support for this comes from an analysis of the sensitivity to initial conditions in the generic game.  Consider, for instance, some particular instantiation of Generic Nim. Now imagine perturbing that instantiation by a single point (i.e., designating one P-position in that instantiation to be an automatic N-position). We can then determine what fraction of P-positions within each loser sheet of the first instantiation will have changed location as a result of the perturbation. Fig. \ref{f13} plots this fractional change as a function of $x$-level. Note that after only about 15 iterations of the recursion operator nearly half of all P-positions in the entire game have moved as a result of this single perturbation.  Thus we see here a form of sensitivity to initial conditions in the generic game. This in turn suggests that long-range correlations between perturbations made at different points will be weak, in further support of the argument for why the pass game exhibits the same overall dynamical behavior as that of the generic game, as was noticed for generic games in \cite{FrL07}.

\section{Concluding remarks}

Our analysis has shown the value of a dynamical-systems approach for understanding combinatorial games.  First we showed that two apparently distinct modifications of games, passes and automatic N-positions (i.e., generic games), are mathematically equivalent in their effect on the recursion operators of the underlying games, providing a critical insight into what has been an open question in combinatorial game theory.  Based on this analysis,  we expect that many other modifications, such as the addition of multiple passes or the introduction of automatic N-positions, will also lead to comparable underlying recursion relations.  In addition, the structural stability (existence of an attracting fixed point for the smoothed system) of the system suggests that all such games should have similar complexity of overall behavior, thus unifying the analysis of a wide variety of game variants.

We have also shown that combinatorial games provide a fascinating set of problems for dynamicists  that can be analyzed by the standard ``tools of the trade,'' but with some new issues and twists (e.g., the computation of lyapunov-like exponents for game operators remains an open question).  Finally, we anticipate that the converse problems may also prove to be of great interest, namely, the use of combinatorial games to understand well known dynamical systems.  

\section{Acknowledgements}

This paper is dedicated to the memory of David Gale, who originally suggested this problem to us.
EJF's research has been supported in part by the NSF under grant CDI-0835706.


\begin{thebibliography}{18}
\bibitem{3Dnote} A three-dimensional position space is not an inherent requirement of this methodology, though it does prove  convenient for visualization purposes.
\end{thebibliography}

\end{document}